\documentclass[12pt]{article}
\usepackage{amssymb,amsfonts,amsmath, psfrag,eepic,colordvi}

\topskip  
\numberwithin{equation}{section}
\newtheorem{theorem}{Theorem}[section]
\newtheorem{example}[theorem]{Example}

\newtheorem{corollary}[theorem]{Corollary}
\newtheorem{definition}[theorem]{Definition}
\newtheorem{conjecture}[theorem]{Conjecture}
\newtheorem{remark}[theorem]{Remark}
\newtheorem{lemma}[theorem]{Lemma}

\begin{document}
\parskip 6pt

\pagenumbering{arabic}
\def\sof{\hfill\rule{2mm}{2mm}}
\def\ls{\leq}
\def\gs{\geq}
\def\SS{\mathcal S}
\def\qq{{\bold q}}
\def\MM{\mathcal M}
\def\TT{\mathcal T}
\def\EE{\mathcal E}
\def\lsp{\mbox{lsp}}
\def\rsp{\mbox{rsp}}
\def\pf{\noindent {\it Proof.} }
\def\mp{\mbox{pyramid}}
\def\mb{\mbox{block}}
\def\mc{\mbox{cross}}
\def\qed{\hfill \rule{4pt}{7pt}}
\def\block{\hfill \rule{5pt}{5pt}}

\begin{center}
{\Large\bf On a refinement of  Wilf-equivalence for  permutations  } \vskip 6mm
\end{center}

\begin{center}
{\small   Huiyun Ge,   Sherry H.F. Yan\footnote{ Corresponding author.
},   Yaqiu Zhang \\[2mm]
 Department of Mathematics, Zhejiang Normal University, Jinhua
321004, P.R. China
\\[2mm]
 huifangyan@hotmail.com
  \\[0pt]
}
\end{center}

\noindent {\bf Abstract.}  Recently, Dokos et al. conjectured that
for all $k, m\geq 1$, the patterns $ 12\ldots k(k+m+1)\ldots
(k+2)(k+1) $ and $(m+1)(m+2)\ldots  (k+m+1)m\ldots 21 $ are
$maj$-Wilf-equivalent. In this paper, we confirm this conjecture for
all $k\geq 1$ and $m=1$. In fact, we construct a descent  set
preserving bijection between $ 12\ldots k (k-1) $-avoiding
permutations and $23\ldots k1$-avoiding permutations for all $k\geq
3$. As a corollary, our bijection  enables  us to  settle a
conjecture of  Gowravaram and Jagadeesan concerning the
Wilf-equivalence for  permutations   with   given descent
sets.

\noindent {\sc Key words}: $maj$-Wilf-equivalent, pattern avoiding permutation, bijection.

\noindent {\sc AMS Mathematical Subject Classifications}: 05A05,
05C30.


\section{Introduction}
 Denote by $\mathcal{S}_n$  the set of all
permutations on $[n]$.  Given a permutation $\pi=\pi_1\pi_2\ldots
\pi_n\in \mathcal{S}_n$ and a permutation $\tau=\tau_1\tau_2\ldots
\tau_k\in \mathcal{S}_k$, we say that $\pi$ contains the {\em pattern}
$\tau$ if there exists a subsequence $\pi_{i_1}\pi_{i_2}\ldots
\pi_{i_k}$ of $\pi$ that is order-isomorphic to $\tau$. Otherwise,
$\pi$ is said to  {\em avoid } the pattern $\tau$ or be {\em
$\tau$-avoiding}. Denote by $\mathcal{S}_n(\tau)$ the set of all $\tau$-avoiding permutations in $\mathcal{S}_n$.
Pattern avoiding permutations have been extensively studied  over last decade.
For a thorough summary of the
current status of research,
see B\'{o}na's  book \cite{bona} and Kitaev's book \cite{kitaev}.

If two patterns $\sigma, \tau\in \mathcal{S}_m$ are said to be {\em Wilf-equivalent}  if and only if $|\mathcal{S}_n(\sigma)|=|\mathcal{S}_n(\tau)|$. A permutation statistic is defined to be a function $s: \mathcal{S}_n\rightarrow T$, where $T$ is any fixed set. The most studied statistics include the inversion number and the major index. Let $\pi=\pi_1\pi_2\ldots \pi_n\in \mathcal{S}_n$. The set of {\em inversions} of $\pi$ is
$$
\mathcal{I}(\pi)=\{(i,j)| i<j\, \mbox{and} \, \pi_i>\pi_j\}.
$$
The inversion number of $\pi$, denoted by $inv(\pi)$,  is   the cardinality of $\mathcal{I}(\pi)$. The {\em decent set} of $\pi$ is
$$
\mathcal{D}(\pi)=\{i| \pi_i>\pi_{i+1}\}.
$$
The {\em ascent set} of $\pi$ is
$$
\mathcal{A}(\pi)=\{i| \pi_i<\pi_{i+1}\}.
$$
The {\em major index} of $\pi$, denoted by $maj(\pi)$, is  given by  $maj(\pi)=\sum_{i\in \mathcal{D}(\pi)}i$.

 Given a permutation statistic $s$,   we say that $\sigma$ and $\tau$ are $s$-wilf-equivalent if there exists a bijection $\Theta: \mathcal{S}_n(\sigma)\rightarrow \mathcal{S}_n(\tau)$ such that $s(\pi)=s(\Theta(\pi))$ for all $\pi\in \mathcal{S}_n(\sigma)$. In other words, the statistic $s$ is equally distributed on the sets $\mathcal{S}_n(\sigma)$ and $\mathcal{S}_n(\tau)$. This refinement of Wilf-equivalence for patterns of length $3$ has been extensively studied, see \cite{Bloom1, Bloom2, Cla, deu, Eli, Rob}. However, little is known about permutation statistics and patterns of length $4$ or greater. Recently, Dokos et al. \cite{dd} posed the following two conjectures on the $maj$-Wilf-equivalence for patterns of length 4 or greater.

 \begin{conjecture}\label{co1}
 ( \cite{dd}, Conjecture 2.6)
 For all $k, m\geq 1$, the patterns $ 12\ldots k(k+m+1) \ldots (k+2)(k+1) $ and $(m+1)(m+2)\ldots  (k+m+1)m\ldots 21 $ are $maj$-Wilf-equivalent.
 \end{conjecture}

 \begin{conjecture}\label{co0}
  ( \cite{dd}, Conjecture 2.8)
 The major index is equally distributed on the sets $\mathcal{S}_n(2413)$, $\mathcal{S}_n(1423)$ and $\mathcal{S}_n(2314)$
 \end{conjecture}

Recently,  Bloom \cite{Bloom3} confirmed Conjecture \ref{co0} by providing   descent set preserving   bijections between the set  $\mathcal{S}_n(2413)$ and the set $\mathcal{S}_n(1423)$,  and between the set $\mathcal{S}_n(2413)$ and the set $\mathcal{S}_n(2314)$.   In their paper \cite{dd},   Dokos et al. showed that Conjecture \ref{co1} is   true for $m=k=1$.  The main purpose of this paper is  to confirm   Conjecture \ref{co1}    for all $k\geq 1$ and $m=1$. Actually, we obtain the following stronger result.
   \begin{theorem}\label{th1}
For $k\geq 3$, there exists a descent set preserving bijection between the set  $\mathcal{S}_n(12\ldots k(k-1))$ and the set $\mathcal{S}_n(23\ldots k 1)$.
 \end{theorem}

Denote by $J_k=12\ldots k$,  $F_k=23\ldots k1$ and  $G_k=12\ldots
k(k-1)$, respectively. Give a permutation
$\pi=\pi_1\pi_2\ldots \pi_n$,  suppose that $\mathcal{D}(\pi)=\{i_1,
i_2, \ldots, i_s\}$. Then we call the subsequence
$\pi_{1}\pi_{2}\ldots \pi_{i_1}$ the first {\em block} of $\pi$, the
subsequence $\pi_{i_1+1}\pi_{i_1+2}\ldots \pi_{i_2}$ the second
block of $\pi$, and so on.
 We say that a permutation $\pi=\pi_1\pi_2\ldots \pi_n$ contains an occurrence of  $H_k$ if  there exists indices $i_1<i_2<\ldots <i_k $   such that the subsequence  $\pi_{i_1} \pi_{i_2} \ldots  \pi_{i_k} $ is isomorphic to $J_k$ and  entries $\pi_{i_{k-1}}$ and $\pi_{i_k}$ belong to two different blocks. That is,   there exists a $j\in  \mathcal{D}(\pi)$ with $i_{k-1}\leq j<i_k$.  Otherwise, we say that $\pi$ avoids $H_k$. For example, the subsequence $1357 9$  of the permutation  $\pi=13576894(10)2(11)\in \mathcal{S}_{11}$  is an occurrence of $H_5$, while the subsequence $13569$ is not an occurrence of $H_5$.
    We say that a permutation $\pi=\pi_1\pi_2\ldots \pi_n$ contains an occurrence of  $Q_k$ if  there exists indices $i_1<i_2<\ldots <i_k $ such that the subsequence $\pi_{i_1}\pi_{i_2}\ldots \pi_{i_k} $ is isomorphic to $G_k$ and $ \pi_{i_{k-1}}<\pi_{i_{k-1}+1}<\ldots< \pi_{i_k-1}>\pi_{i_k}$. Otherwise, we say that  $\pi$ avoids $Q_k$. For example, the subsequence $13586$  of the permutation  $\pi=1358(10)67492(11)\in \mathcal{S}_{11}$  is an occurrence of $Q_5$, while the subsequence $13587$ is not an occurrence of $Q_5$.

 In order to prove Theorem \ref{th1}, we obtain the following two theorems.
\begin{theorem}\label{th2}
For $k\geq 3$, there is a bijection $f$ between the set $\mathcal{S}_n(G_k)$ and   the set $\mathcal{S}_n(H_k, Q_k)$ such that for any $\pi\in \mathcal{S}_n(G_k)$, we have $\mathcal{D}(\pi)=\mathcal{D}(f(\pi)) $.
\end{theorem}

\begin{theorem}\label{th3}
For $k\geq 3$, there is a bijection $\Phi$ between the set $\mathcal{S}_n(F_k)$  and the set  $\mathcal{S}_n(H_k, Q_k)$ such that for any $\pi\in \mathcal{S}_n(F_k)$, we have $\mathcal{D}(\pi)=\mathcal{D} (\Phi(\pi))$.
\end{theorem}
Combining Theorems \ref{th2} and \ref{th3}, we are led to Theorem \ref{th1}.

Given a positive  integer $t$, Let $D^t_{n}=\{i| 1\leq i\leq n-1 \, \mbox{and} \,i\equiv 0\, \mbox{mod} \,t\}$. Denote by $\mathcal{S}^t_n(12\ldots k(k-1))$ (resp. $\mathcal{S}^t_n(23\ldots k 1)$) the set of permutations  $\pi\in \mathcal{S}_n(12\ldots k(k-1))$ (resp. $\pi\in \mathcal{S}^t_n(23\ldots k 1)$) with $\mathcal{D}(\pi)=D^t_{n}$.    From Theorem \ref{th1}, we obtain the following result as   conjectured   by Gowravaram and   Jagadeesan \cite{gj}.
\begin{corollary} (\cite{gj}, Conjecture 6.2)
 For $t\geq 1$ and $k\geq 3$, we have $|\mathcal{S}^t_n(12\ldots k(k-1))|=|\mathcal{S}^t_n(23\ldots k 1)|$.
\end{corollary}

\section{Proof of Theorem \ref{th2}}

We begin with some definitions and notations.
 An entry of a permutation is said to have {\em rank }$k$ if  the length
  of the longest increasing subsequence that ends in that entry is  $k$.
 We now construct a map  $f$ from the set $\mathcal{S}_n(G_k)$ to the set $\mathcal{S}_n(H_k, Q_k)$.  The map $f$ is a slight modification  of a classic bijection, which is given by West \cite{west} to prove the equality $|\mathcal{S}_n(J_k)|=|\mathcal{S}_n(G_k)|$ for all $k\geq 3$. Recently, Bona \cite{bona2} proved that West's bijection also induces a bijection between $12\ldots k$-avoiding alternating permutations and $12\ldots k(k-1)$-avoiding alternating permutations, thereby proving  generalized versions of some conjectures of Lewis \cite{Lewis2}.

 Let $\pi\in \mathcal{S}_n(G_k)$.  In order to obtain $f(\pi)$, we leave all entries of $\pi$ that are of rank $k-2$ or less in their place and rearrange the entries of rank    $k-1$ or higher. Let $B_1, B_2, \ldots, B_s$ be the blocks of $\pi$ that are listed from left to right. Let $P_i$ be the set of positions of $\pi$ in which,  an entry that has  rank   $k-1$ or higher and  belongs to  the block $B_i$,   is located,  and let $R$ be the set of entries of $\pi$ that are of rank   $k-1$ or higher.  Now we fill the positions of $P_i's$ as follows.
\begin{itemize}
\item[Step 1.] Choose $|P_1|$ largest   entries from $R$ that are larger than the closest entry of rank $k-2$ on the left of the positions of $P_1$,   and   fill the positions of $P_1$ with the selected   entries  from left to right in increasing order.
    \item[Step 2.]
  Choose $|P_2|$ largest   entries from $R$ that  have not been placed yet and are larger than the closest entry of rank $k-2$ on the left of the positions of $P_2$. Fill the positions of $P_2$ with the selected entries from left to right in increasing order.
 \item[Step 3.]    Fill the positions of  $P_{3}$, $P_{4}$, $\ldots$, $P_s$ as in Step 2.

   \end{itemize}
Let $f(\pi)$ be the obtained permutation.
\begin{example}

 Consider  $\pi=13576894(10)2(11)\in \mathcal{S}_{11}(G_6)$. Then we have $B_1= 1 3 5 7 $, $B_2= 6 8 9 $, $B_3= 4 (10) $ and $B_4= 2(11) $. Moreover, we have $P_1=\emptyset$, $P_2=\{6,7\}$, $P_3=\{9\}$,  $P_4=\{11\}$ and $R=\{8,9,10,11\}$. According to the definition of $f$, we have
 $f(\pi)=13576(10)(11)4928$.
\end{example}

Since the existence of $\pi$ shows that there is at least one way to assign the entries of $R$ to the positions of $P_i$,  the definition of $f$ always enables us to create   $f(\pi)$.

 Notice that if entry $\pi_i$ of $\pi$ has rank $k-2$ or less, then $\pi_i$ do not move in the above procedure, and the rank of $\pi_i$ do not change. If   entry $\pi_i$ of $\pi$ has rank $k-1$ or higher, then $\pi_i$ may have moved and the rank of $\pi_i$ in $f(\pi)$ is $k-1$ or higher.   We claim that if $\pi_{i-1}>\pi_i$, then the rank of $\pi_i$ is $k-2$ or less. If not,   the longest increasing subsequence that ends in $\pi_i$ combining with $\pi_{i-1}$ would  form a $G_k$ in $\pi$. This contradicts the fact that $\pi$ avoids $G_k$.

 Now we proceed to show that $\mathcal{D} (\pi)=\mathcal{D} (f(\pi))$.    Let $f(\pi)=\sigma_1\sigma_2\ldots \sigma_n$.  If $\pi_{i-1}>\pi_{i}$, then the rank of $\pi_i$ is $k-2$ or less and do not move. This implies that $\pi_{i}=\sigma_i$ and $\sigma_i$ has rank $k-2$ or less. If $\pi_{i-1}$ is of rank $k-2$ or less, then we have $\sigma_{i-1}=\pi_{i-1}$. In this case, we have $\sigma_{i-1}=\pi_{i-1}>\pi_i=\sigma_i$. If  $\pi_{i-1}$ is of rank $k-1$ or higher, then $\sigma_{i-1}$ is of rank $k-1$ or higher. Since $\sigma_{i}$ is of rank $k-2$ or less,   we have  $\sigma_{i-1}>\sigma_i$. Thus, we have concluded that if  $\pi_{i-1}>\pi_{i}$, then $\sigma_{i-1}>\sigma_i$.

Next we aim to show that if $\pi_{i-1}<\pi_{i}$, then we have $\sigma_{i-1}<\sigma_i$. We have three cases. If $\pi_i$ is of rank $k-2$ or less in $\pi$, then then the rank of $\pi_{i-1}$ is also $k-2$ or less. In this case, we have $\sigma_{i-1}=\pi_{i-1}<\pi_i=\sigma_i$. If both $\pi_i$  and $\pi_{i-1}$ are of rank $k-1$ or higher, then according to the definition of $f$, we have $\sigma_{i-1}<\sigma_{i}$ . If $\pi_i$ has rank $k-1$ or higher and  $\pi_{i-1}$ is of rank $k-2$ or less, then the rank of $\sigma_{i-1}$ is $k-2$ or less and $\sigma_i$ is of rank $k-1$ or higher. This implies that $\sigma_{i-1}<\sigma_i$. Thus, we have concluded that if  $\pi_{i-1}<\pi_{i}$, then $\sigma_{i-1}<\sigma_i$.  Therefore,  we have   $\mathcal{D}(\pi)=\mathcal{D}(f(\pi))$.

   Notice that $f(\pi)$ avoids $H_k$ since the existence of such a pattern in $f(\pi)$ would mean that the last two entries of that pattern were not placed  according to the rule specified above.  Moreover,  we have that $f(\pi)$ avoids $Q_k$. If not, suppose that $\sigma_{i_1}\sigma_{i_2}\ldots \sigma_{i_k}$ is such a $Q_k$. Then we have $\sigma_{i_k-1}>\sigma_{i_k}$ and $\sigma_{i_k}$ has rank $k-1$ or higher. Since $\mathcal{D}(\pi)=\mathcal{D}(f(\pi))$, we have $\pi_{i_k-1}>\pi_{i_k}$. Recall that if $\pi_{i-1}>\pi_i$, then both $\pi_i$ and $\sigma_i$ have rank  $k-2$ or less. This implies that $\sigma_{i_k}$ has rank $k-2$ or less, which  contradicts the fact that $\sigma_{i_k}$ has rank $k-1$ or higher. Thus, we deduce that $f(\pi)$ avoids $Q_k$.

In order to show that the map $f$ is a bijection, we  construct a
map  $g$ from the set $\mathcal{S}_n(H_k, Q_k)$ to the set
$\mathcal{S}_n(G_k)$. Let $\sigma=\sigma_1\sigma_2\ldots\sigma_n\in
\mathcal{S}_n(H_k, Q_k)$. We aim to obtain $g(\sigma)$ by leaving
all entries of $\sigma$ that are of rank $k-2$ or less in their
place and rearranging the entries of rank    $k-1$ or higher. Label
the blocks of $\sigma $ from left to right by  $B_1, B_2, \ldots,
B_s$.  Let $P_i$ be the set of positions of $\pi$ in which an entry,
that has  rank  $k-1$ or higher  and belongs to the block $B_i$,  is
located,  and let $R$ be the set of entries of $\pi$ that are of
rank   $k-1$ or higher.  Now we fill the positions of $P_i$ as
follows.
\begin{itemize}
\item[Step 1.] Choose $|P_1|$ smallest  entries from $R$ that are larger than the closest entry of rank $k-2$  on the left of the positions of $P_1$,   and   fill the positions of $P_1$ with the selected   entries  from left to right in increasing order.
    \item[Step 2.]
  Choose $|P_{2}|$ smallest   entries from $R$ that  have not been placed yet and are larger than the closest entry of rank $k-2$ on  the left of the positions of $P_2$. Fill the positions of $P_2$ with the selected entries from left to right in increasing order.
 \item[Step 3.]    Fill the positions of  $P_{3}$,  $P_4$,  $\ldots$, $P_s$ as in Step 2.

   \end{itemize}
Let $g(\sigma)$ be the obtained permutation.
\begin{example}
 Consider  $\pi=13576(10)(11)4928\in \mathcal{S}_{11}(H_6, Q_6)$. Then we have $B_1= 1 3 5 7 $,
 $B_2= 6 (10)(11) $, $B_3= 4 9$ and $B_4= 2 8 $. Moreover, we have $P_1=\emptyset$, $P_2=\{6,7\}$, $P_3=\{9\}$,
  $P_4=\{11\}$ and $R=\{8,9,10,11\}$. According to the definition of $g$, we have
 $g(\pi)=13576894(10)2(11)$.
\end{example}

 Since the existence of $\sigma$ shows that there is at least one way to assign the entries of $R$ to the positions of $P_i$,  the definition of $g$ always enables us to create   $g(\sigma)$.

 Notice  that if entry $\sigma_i$ has rank $k-2$ or less, then $\sigma_i$ does not move in the above procedure, and the rank of $\sigma_i$ do not change. If   entry $\sigma_i$ has rank $k-1$ or higher, then $\sigma_i$ may have moved and the rank of $\sigma_i$ in $g(\sigma)$ is $k-1$ or higher.   We claim that if $\sigma_{i-1}>\sigma_i$, then the rank of $\sigma_i$ is $k-2$ or less. If not,    the longest increasing subsequence that ends in $\sigma_i$ combining with $\sigma_{i-1}$ would  form a $Q_k$ in $\sigma$.

 By the similar reasoning as in the proof of the equality $\mathcal{D}(\pi)=\mathcal{D}(f(\pi))$,
 one can verify that $\mathcal{D}(\sigma)=\mathcal{D}(g(\sigma))$.  Now we proceed to show that $g(\sigma)$ avoids $G_k$.
 Let $g(\sigma)=\pi_1\pi_2\ldots \pi_n$.  Suppose that the the  subsequence  $\pi_{i_1}\pi_{i_2}\ldots \pi_{i_k}$ is a pattern $G_k$
 in $\pi$ with $i_1<i_2<\ldots<i_k$.  Without loss of generality, assume that $\pi_{i_{k-2}}$ has rank $k-2$.  Clearly, both $\pi_{i_{k-1}}$
 and $\pi_{i_k}$ have rank $k-1$ or higher in $\pi$.   Suppose that
 $i_{k-1}\in P_j$ for some $j$, and $\sigma_s$ is the closest entry of rank $k-2$  on the left of the positions of
 $P_j$ in $\sigma$. Recall that the map $g$ does not change the position of  entry $\sigma_i$  that has  rank $k-2$ or
 less, and the rank of $\sigma_i$ does not change in $\pi$. So we
 have $\pi_s=\sigma_{s}$ and $\pi_{i_{k-2}}=\sigma_{i_{k-2}}$, and the rank of $\pi_s$ (resp. $\sigma_{i_{k-2}}$) is $k-2$ in $\pi$ (resp. $\sigma$).  Moreover, since  $\sigma_s$ is the closest entry of rank $k-2$  on the left of the positions of
 $P_j$ in $\sigma$,   we have $i_{k-2}\leq s$. This
 implies that $\sigma_s \leq  \sigma_{i_{k-2}}=\pi_{i_{k-2}}$.  Then, we have
 $\pi_{i_{k}}>\sigma_s$, which contradicts the selection of
 $\pi_{i_{k-1}}$ when filling the positions of $P_j$. Hence, we have
 $g(\sigma)\in \mathcal{S}_n(G_k)$.

 In order to show that the map $f$ is a bijection, it suffices to show that the maps $f$ and $g$ are inverses of each other.
  First, we wish to prove that for any $\pi\in \mathcal{S}_n(G_k)$, we have
  $g(f(\pi))=\pi$.  Since $\mathcal{D}(f(\pi))=\mathcal{D}(\pi)$,  $\pi$ and  $f(\pi)$ have the same number of blocks.  Suppose that $B'_1, B'_2, \ldots, B'_s$ are the blocks of $f(\pi)$, that are listed from
   left to right. Let $P'_i$ be the set of positions of $f(\pi)$ in which an entry that has rank  $k-1$ or higher and belongs to
   the block $B'_i$, is located,  and let $R'$ be the set of entries of $f(\pi)$ that are of rank   $k-1$ or higher.   Recall that if the entry $\pi_i$ of $\pi$ has rank $k-2$ or less, then the map $f$ does not
    change the position of  $\pi_i$, and the rank of $\pi_i$ do not change. If   entry $\pi_i$ of $\pi$ has rank $k-1$ or higher,
    then $\pi_i$ may have moved and the rank of $\pi_i$ in $f(\pi)$ is $k-1$ or higher.   So we have $P_i=P'_i$ and $R=R'$.
  Since $\pi$ avoids $G_k$, the positions of $P_1$ in $\pi$ are filled with $|P_1|$ smallest elements of $R$
   in increasing order which are larger than the closet entry of rank $k-2$ on the left of the positions of $P_1$.
   The positions of $P_2$ are filled with the next $|P_2|$ smallest entries of $R$ in increasing order that have not been placed and
    larger than the closet entry of rank $k-2$  on the left of the positions of $P_2$. And the positions of $P_3, \ldots, P_s$ are filled
    in the same manner as the positions of $P_2$. Thus, according to the definition of $g$, it is easy to check that $g(f(\pi))=\pi$.

Our next goal is to   show that for any $\sigma\in
\mathcal{S}_n(H_k, Q_k)$, we have $f(g(\sigma))=\sigma$.   Since
$\mathcal{D}(g(\sigma))=\mathcal{D}(\sigma)$, $\sigma$ and
$f(\sigma)$ have the same number of blocks.  Suppose that $B'_1,
B'_2, \ldots, B'_s$ are the blocks of $f(\sigma)$, that are listed
from
   left to right.
Let $P'_i$ be the set of positions of $f(\pi)$ in which an entry
that has rank  $k-1$ or higher and belongs to the block $B'_i$, is
located,  and let $R'$ be the set of entries of $f(\pi)$ that are of
rank   $k-1$ or higher.  Recall that if entry $\sigma_i$ of $\sigma$
has rank $k-2$ or less, then the map $g$ does not change the
position of  $\sigma_i$,  and the rank of $\sigma_i$ do not change.
If   entry $\sigma_i$ of $\sigma$ has rank $k-1$ or higher, then
$\sigma_i$ may have moved and the rank of $\sigma_i$ in $g(\sigma)$
is $k-1$ or higher.  So we have $P_i=P'_i$ and $R=R'$.
  Since $\sigma$ avoids $H_k$, the positions of $P_1$ in $\sigma$ are  filled with $|P_1|$
  largest elements of $R$   in increasing order which are larger than the closet entry of rank $k-2$
   on the left of the positions of $P_1$. The positions of $P_2$ are filled with the next $|P_2|$ largest
    entries of $R$ in increasing order that have not been placed and   larger than the closet entry of rank $k-2$ on the left of  the positions of $P_2$. And the positions of $P_3, \ldots, P_s$ are filled in the same manner as the positions of $P_2$. Thus, according to the definition of $f$, it is easy to check that $f(g(\sigma))=\sigma$.

  \section{Proof of Theorem \ref{th3}}
  Let us begin with some necessary definitions and notations.  We draw Young
diagrams in English notation,  and   number columns from left to right and  rows from bottom to up.
  For example,   the square $(1, 2)$ is the second
square in the bottom row of a Young diagram.

A {\em  transversal} of  a Young diagram $\lambda=(\lambda_1\geq \lambda_2\geq \ldots\geq \lambda_n)$  is a filling  of the squares of $\lambda$ with $1's$ and $0's$ such that every row and column contains  exactly one $1$.    Denote by $T=\{(t_i,i)\}_{i=1}^{n}$ the transversal in which the square $(t_i, i)$ is filled with a $1$ for all $i\leq n$. For example the transversal $T=\{(1,1),  (2,3),  (3,2), (4, 4), (5,5)\}$ of a Young diagram $(5,4,3,3,1)$ is illustrated as Figure \ref{fig1}.

  \begin{figure}[h,t]
\begin{center}
\begin{picture}(50,30)
\setlength{\unitlength}{1.5mm}

\put(0,0){\framebox(5,5){$1$} }

 \put(0,5){\framebox(5,5){$0$}  } \put(5,5){\framebox(5,5)  {$0$}}  \put(10,5){\framebox(5,5) {$1$} }
 \put(0,10){\framebox(5,5){$0$}  } \put(5,10){\framebox(5,5)  {$1$}}  \put(10,10){\framebox(5,5) {$0$} }

 \put(0,15){\framebox(5,5){$0$}  } \put(5,15){\framebox(5,5)  {$0$}}  \put(10,15){\framebox(5,5) {$0$} }\put(15,15){\framebox(5,5) {$1$} }

     \put(0,20){\framebox(5,5){$0$}  } \put(5,20){\framebox(5,5)  {$0$}}  \put(10,20){\framebox(5,5) {$0$} }\put(15,20){\framebox(5,5) {$0$} }

     \put(20,20){\framebox(5,5){$1$}  }

  \end{picture}
\end{center}
\caption{ The transversal $T=\{(1,1),  (2,3),  (3,2), (4, 4), (5,5)\}$.}\label{fig1}
\end{figure}

In this section, we will  consider permutations as permutation matrices. Given a permutation $\pi=\pi_1\pi_2\ldots\pi_n\in \mathcal{S}_n$, its corresponding {\em permutation matrix} is a transversal of the square shape $\lambda_1=\lambda_2=\ldots=\lambda_n=n$ in which the square   $(\pi_i,i)$ is filled with  a $1$ for all $1\leq i\leq n$ and all the other squares  are filled with $0's$.

The   notion of pattern  avoidance is extended to transversal of a Young diagram in  \cite{BWX}.
Given a permutation $\alpha$ of $\mathcal{S}_m$,  let   $M$ be its permutation matrix. A transversal $L$ of a Young diagram $\lambda$  will be said to contain   $\alpha$ if there exists two subsets of the index set $[n]$, namely,  $R=\{r_1<r_2<\ldots<r_m\}$ and $C=\{c_1< c_2< \ldots< c_m\}$, such that  the   matrix on $R$ and $C$ is a copy of $M$ and each of the squares $(r_j,c_j)$ falls within the Young diagram.

\subsection{The map $\Phi$   from the set $\mathcal{S}_n(F_k)$ to the set $\mathcal{S}_n(H_k, Q_k)$}

 Before we describe  the map $\Phi$, let us review a  transformation   $\theta$ introduced in \cite{BWX}

 Let $\pi=\{(\pi_1, 1), (\pi_2, 2), \ldots, (\pi_n,n)\}$.
 Suppose that  $G$ is the submatrix  of $\pi$ at columns $c_1<c_2<\ldots<c_{k-1}<c_k$ and rows $r_1< r_2< \ldots< r_{k-1}< r_k$, which is isomorphic to $J_k$.
 In other words,
 the square  $(r_i, c_i)$ is filled with $1$ for all $i=1, 2, \ldots, k$. Let $\theta(G)$ be the submatrix at the same rows and columns as $G$,  such that the squares
    $(r_2, c_1)$, $(r_3, c_2)$, $\ldots$, $(r_{k}, c_{k-1})$, $(r_1, c_k)$ are filled with $1's$ and all the other squares are filled with $0's$. Clearly, $\theta(G)$ is isomorphic to $F_k$.

Based on the transformation $\theta$, we define the following three transformations, which will play an essential role in the construction of the map $\Phi$.

 Suppose that  $G$ is the submatrix  of $\pi$ at columns $c_1<c_2<\ldots<c_{k-1}<s<s+1<\ldots<c_{k}-1<c_k$ and rows $r_1< r_2< \ldots< r_{k-1}<r_k> \pi_{c_k-1}> \ldots> \pi_{s+1}>\pi_s$, in which the squares $(r_i, c_i)$ are filled with $1's$ for all $i=1, 2, \ldots, k$.
 Let $\alpha(G)$ be the submatrix at the same rows and columns as $G$,  such that the squares
    $(r_2, c_1)$, $(r_3, c_2)$, $\ldots$, $(r_{k}, c_{k-1})$, $(r_1, s)$, $(\pi_{s}, s+1)$, $\ldots$,
    $ (\pi_{c_k-1}, c_k)$ are filled with $1's$ and all the other squares are filled with $0's$.
    Clearly,    the submatrix at columns $c_1<c_2<\ldots<c_{k-1}<s$ and rows $r_1<r_2<\ldots<r_{k-1}<r_k$ is isomorphic to $F_k$.

Suppose that $G$ is the submatrix of $\pi$ at  columns $c_1<c_2<\ldots<c_{k-1}< c_k<c_k+1<\ldots<t-1<t$ and rows $r_1< r_2<\ldots<r_{k}>\pi_{c_k+1}> \ldots>\pi_{t-1}> \pi_{t}$,
 in which the squares $(r_i, c_i)$ are filled with $1's$ for all $i=1, 2, \ldots, k$.   Define $\beta(G)$ to be the submatrix at the same  columns  and  rows as $G$,  such that the squares
    $(r_2, c_1)$, $(r_3, c_2)$, $\ldots$, $(r_{k}, c_{k-1})$, $(\pi_{c_k+1}, c_k)$,   $\ldots$, $(\pi_{t}, t-1)$, $(r_1, t)$
are filled with $1's$ and all the other squares are filled with $0's$. Clearly, the submatrix at columns $c_1<c_2<\ldots<c_{k-1}<t$ and rows $r_1<r_2<\ldots<r_{k-1}<r_k$ is isomorphic to $F_k$.

 {\noindent \em The transformation $\phi$:  }
Suppose that $\pi $ is a permutation in $\mathcal{S}_n$.  First, find the highest square  $ (p_1, q_1)$ containing a $1$,
 such that there is an $H_k$ or $Q_k$ in $\pi$ in which the $1$ positioned at the square$(p_1, q_1)$ is the leftmost entry.
 Then, find the leftmost square  $(p_2, q_2)$ containing a $1$,  such that  there is an $H_k$ or $Q_k$ in $\pi$ in which
  the $1's$ positioned at the squares $(p_1, q_1)$   and $(p_2, q_2)$  are the leftmost two $1's$.
  Find $(p_3, q_3), (p_4, q_4), \ldots, (p_{k-1}, q_{k-1})$ one by one as $(p_2, q_2)$.

 If there is an $H_k$ in which  the $1's$ positioned at the squares $(p_1, q_1),   (p_2, q_2),   \ldots $ $  (p_{k-1}, q_{k-1})$ are the
 leftmost $k-1$ $1's$, then find the highest square   $(p_k, q_k)$  containing a $1$, such that  the $1's$ positioned at
  the squares   $(p_1, q_1),   (p_2, q_2),   \ldots $ $  (p_k, q_k)$  form an  $H_k$.
   Find  the largest $s$ such that $s-1\in \mathcal{D}(\pi)$ and $q_{k-1}<s< q_{k}$.
    Now we proceed to   construct a permutation $\phi(\pi)$  by the following procedure.
    \begin{enumerate}
    \item[Case 1.] $q_k=n$ or $\pi_{q_k-1}>\pi_{q_k+1}$.    Let $G$ be the submatrix  of $\pi$ at columns
    $q_1<q_2<\ldots<q_{k-1}<s<s+1<\ldots<q_k-1<q_k$ and rows $p_1< p_2< \ldots <p_k> \pi_{q_k-1}>\ldots> \pi_{s+1}> \pi_{s}$.  Replace $G$ by $\alpha(G)$ and leave all the other squares fixed.

    \item[Case 2.] $\pi_{q_k-1}<\pi_{q_k+1}$.    Find the least  $t$ such that $t> q_k$ and $t\in \mathcal{A}(\pi)$. If such $t$ does not exist, set $t=n$. Let $G$ be the submatrix of $\pi$ at  columns $q_1<q_2<\ldots<q_{k-1}< q_k<q_k+1<\ldots<t-1<t$ and rows $p_1< p_2<\ldots<p_{k}>\pi_{q_k+1}> \ldots>\pi_{t-1}> \pi_{t}$.
          Replace $G$ by $\beta(G)$  and leave all the other squares fixed.
           \end{enumerate}
   If such an $H_k$ does not exist,  then find the leftmost square  $(p_k, q_k)$ containing a $1$, such that the $1's$ positioned at the squares    $(p_1, q_1), (p_2, q_2) \ldots,   (p_k, q_k)$ form a  $Q_k$.   Construct a permutation $\phi(\pi)$   by the following procedure.

   \begin{enumerate}
   \item[Case 3.] $q_k\in \mathcal{A}(\pi)$.  Let $G$ be the submatrix of $\pi$ at  columns $q_1<q_2<\ldots<q_{k-2}< q_k$ and rows $p_1< p_2<\ldots<p_{k-2}<p_k$.
          Replace $G$ by $\theta(G)$  and leave all the other squares fixed.

            \item[Case 4.] Otherwise, find the least  $t$ such that $t> q_k$ and $t\in \mathcal{A}(\pi)$. If such $t$ does not exist, set $t=n$. Let $G$ be the submatrix of $\pi$ at  columns $q_1<q_2<\ldots<q_{k-2}<q_k< q_k+1<q_{k}+2<\ldots< t-1<t$ and rows $p_1< p_2< \ldots< p_{k-2}< p_k> \pi_{q_k+1}>\pi_{q_k+2}>\ldots>\pi_{t-1}>\pi_t$.  Replace $G$ by $\beta(G)$  and leave all the other squares fixed.
\end{enumerate}

\begin{remark}
Notice that  the definition of $H_k$ ensures that there exists  an
$s$  such that $s-1\in \mathcal{D}(\pi)$ and $q_{k-1}<s\leq q_{k}$.
In fact, we have     $q_{k-1}<s<q_k$.    If not, then the $1's$
positioned at the squares    $(p_2, q_2), (p_3, q_3) \ldots(p_{k-1},
q_{k-1}), (\pi_{q_k-1}, q_k-1), (p_k, q_k)$ would form a $Q_k$,
which contradicts the selection of $(p_1, q_1)$.
\end{remark}

\begin{remark}
We denote the resulting permutation  in Case 1, Case 2, Case 3 and Case 4 by  $\phi_1(\pi), \phi_2(\pi)$,   $\phi_3(\pi)$ and  $\phi_4(\pi)$, respectively.
\end{remark}

It is obvious that the transformation $\phi$   changes every occurrence of $H_k$ (or  $Q_k$) to an occurrence of $F_k$.
 Denote by    $\Phi$   the iterated transformation, that recursively transforms every occurrence of $H_k$ (or  $Q_k$) into $F_k$.

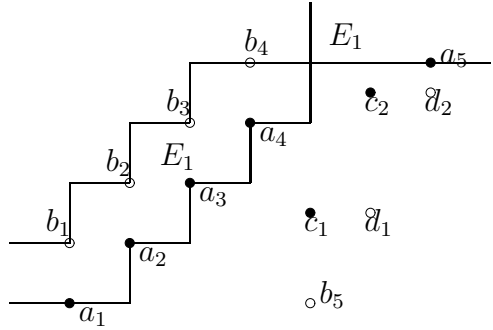
\begin{figure}[h,t]
\begin{center}
\begin{picture}(50,30)
\setlength{\unitlength}{8mm} \linethickness{0.4pt}
\put(0,0){\line(1,0){2}}\put(1,0){\circle*{0.15}}
\put(2,0){\line(0,1){1}}\put(2,1){\circle*{0.15}}\put(2,1){\line(1,0){1}}
\put(3,1){\line(0,1){1}}\put(3,2){\circle*{0.15}}\put(3,2){\line(1,0){1}}
\put(4,2){\line(0,1){1}}\put(4,3){\circle*{0.15}}\put(4,3){\line(1,0){1}}
\put(5,3){\line(0,1){2}}\put(7,4){\circle*{0.15}}
\put(0,1){\line(1,0){1}}\put(1,1){\circle{0.15}}
\put(1,1){\line(0,1){1}}\put(1,2){\line(1,0){1}}\put(2,2){\circle{0.15}}
\put(2,2){\line(0,1){1}}\put(2,3){\line(1,0){1}}\put(3,3){\circle{0.15}}
\put(3,3){\line(0,1){1}}\put(3,4){\line(1,0){1}}\put(4,4){\circle{0.15}}\put(4,4){\line(1,0){4}}
\put(5,0){\circle{0.15}}\put(5,1.5){\circle*{0.15}}\put(6,1.5){\circle{0.15}}
\put(6,3.5){\circle*{0.15}}\put(7,3.5){\circle{0.15}}\put(2.5,2.3){$E_{1}$}\put(5.3,4.3){$E_{1}$}
\put(1.15,-0.3){$a_{1}$} \put(2.15,0.7){$a_{2}$}\put(3.15,1.7){$a_{3}$}\put(4.15,2.7){$a_{4}$}\put(7.15,4){$a_{5}$}
\put(0.6,1.2){$b_{1}$}\put(1.6,2.2){$b_{2}$}\put(2.6,3.2){$b_{3}$}\put(3.9,4.2){$b_{4}$}\put(5.15,0){$b_{5}$}
\put(4.9,1.2){$c_{1}$}\put(5.9,1.2){$d_{1}$}\put(5.9,3.2){$c_{2}$}\put(6.9,3.2){$d_{2}$}
\end{picture}
\end{center}
\caption{The labelling of squares in $\phi_1(\pi)$ for $k=5$.} \label{case1}
\end{figure}

Using the notation of the algorithm for $\phi_1$, we label the squares containing $1's$ in $G$  by $a_1, a_2, \ldots, a_{k-1}, c_1, c_2, \ldots, c_{q_k-s}, a_{k}$,    and   the squares   containing   $1's$  in  $\alpha(G)$  by $b_1, b_2, \ldots, b_k, d_1, d_2, \ldots, d_{q_k-s}$,  from left to right, see Figure \ref{case1} for example.

\begin{figure}[h,t]
\begin{center}
\begin{picture}(50,40)
\setlength{\unitlength}{7mm} \linethickness{0.4pt}
\put(0,0){\line(1,0){2}}\put(1,0){\circle*{0.15}}
\put(2,0){\line(0,1){1}}\put(2,1){\circle*{0.15}}\put(2,1){\line(1,0){1}}
\put(3,1){\line(0,1){1}}\put(3,2){\circle*{0.15}}\put(3,2){\line(1,0){1}}
\put(4,2){\line(0,1){1}}\put(4,3){\circle*{0.15}}\put(4,3){\line(1,0){1}}
\put(5,3){\line(0,1){3}}\put(6,5){\circle*{0.15}}
\put(0,1){\line(1,0){1}}\put(1,1){\circle{0.15}}
\put(1,1){\line(0,1){1}}\put(1,2){\line(1,0){1}}\put(2,2){\circle{0.15}}
\put(2,2){\line(0,1){1}}\put(2,3){\line(1,0){1}}\put(3,3){\circle{0.15}}
\put(3,3){\line(0,1){2}}\put(3,5){\line(1,0){1}}\put(4,5){\circle{0.15}}\put(4,5){\line(1,0){3}}
\put(8,0){\circle{0.15}}\put(5,1.5){\circle*{0.15}}
\put(6,4.2){\circle{0.15}}\put(7,4.2){\circle*{0.15}}\put(4,4){$E_{1}$}\put(5.3,5.3){$E_{1}$}
\put(7,2.5){\circle{0.15}}\put(8,2.5){\circle*{0.15}}
\put(1.15,-0.3){$a_{1}$} \put(2.15,0.7){$a_{2}$}\put(3.15,1.7){$a_{3}$}\put(4.15,2.7){$a_{4}$}\put(6.15,5){$a_{5}$}
\put(0.6,1.2){$b_{1}$}\put(1.6,2.2){$b_{2}$}\put(2.6,3.2){$b_{3}$}\put(3.9,5.2){$b_{4}$}\put(8.15,0){$b_{5}$}
\put(5.15,1.5){$c_{1}$}\put(5.9,3.9){$f_{1}$}\put(6.9,2.2){$f_{2}$}\put(7.1,3.9){$e_{1}$}\put(8.1,2.2){$e_{2}$}

\end{picture}
\end{center}
\caption{The labelling of squares in $\phi_2(\pi)$ for $k=5$.} \label{case2}
\end{figure}
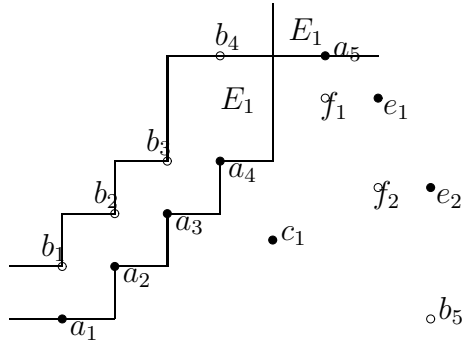

Using the notation of the algorithm for $\phi_2$,  we label the squares containing $1's$ in $G$   by $a_1, a_2, \ldots, a_{k-1},     a_{k}, e_1, e_2, \ldots, e_{t-q_k}$,   and   the squares   containing   $1's$ in  $\beta(G)$  by $b_1, b_2, \ldots, b_{k-1}, f_1, f_2, \ldots, f_{t-q_k}, b_k$, from left to right. We also label the square $(\pi_s, s)$ by $c_1$,  see Figure \ref{case2} for example.

\begin{figure}[h,t]
\begin{center}
\begin{picture}(50,40)
\setlength{\unitlength}{7mm} \linethickness{0.4pt}
\put(0,0){\line(1,0){2}}\put(1,0){\circle*{0.15}}
\put(2,0){\line(0,1){1}}\put(2,1){\circle*{0.15}}\put(2,1){\line(1,0){1}}
\put(3,1){\line(0,1){1}}\put(3,2){\circle*{0.15}}\put(3,2){\line(1,0){1}}
\put(4,2){\line(0,1){1}}\put(4,3){\circle*{0.15}}\put(4,3){\line(1,0){3}}
\put(0,1){\line(1,0){1}}\put(1,1){\circle{0.15}}
\put(1,1){\line(0,1){1}}\put(1,2){\line(1,0){1}}\put(2,2){\circle{0.15}}
\put(2,2){\line(0,1){1}}\put(2,3){\line(1,0){1}}\put(3,3){\circle{0.15}}
\put(3,3){\line(0,1){1}}\put(3,4){\line(1,0){1}}\put(4,4){\circle{0.15}}
\put(4,4){\line(0,1){1}}\put(4,5){\line(1,0){5}}\put(6,5){\circle*{0.15}}\put(7,3){\line(0,1){3}}
\put(7,4){\circle*{0.15}}\put(7,0){\circle{0.15}}\put(5,4){$E_{2}$}\put(7.3,5.3){$E_{2}$}
\put(1.15,-0.3){$a_{1}$} \put(2.15,0.7){$a_{2}$}\put(3.15,1.7){$a_{3}$}\put(4.15,2.7){$a_{4}$}\put(5.6,5.2){$a_{5}$}\put(7.1,3.9){$a_{6}$}
\put(0.6,1.2){$b_{1}$}\put(1.6,2.2){$b_{2}$}\put(2.6,3.2){$b_{3}$}\put(3.6,4.2){$b_{4}$}\put(7.15,0){$b_{6}$}\put(6,5.2){$(b_{5})$}

\end{picture}
\end{center}
\caption{The labelling of squares in $\phi_3(\pi)$ for $k=6$.} \label{case3}
\end{figure}
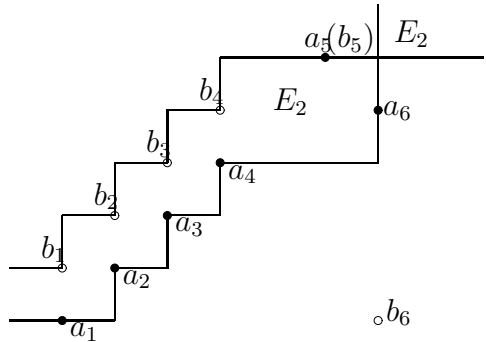

Using the notation of the algorithm for $\phi_3$, we label the squares  containing $1's$ in $G$  by $a_1, a_2, \ldots, a_{k-2}, a_k,$ and
 the squares  containing $1's$ in $\theta(G)$ by $b_1, b_2, \ldots, b_{k-2},   b_k$,   from left to right. We also label the square $(p_{k-1}, q_{k-1})$ by $b_{k-1}$ (or $a_{k-1}$), see Figure \ref{case3} for example.

Using the notation of the algorithm for $\phi_4$, we label the squares  containing $1's$ in $G$  by $a_1, a_2, \ldots, a_{k-2}, a_k,$ $ e_1, e_2, \ldots, e_{t-q_k} $, and  the squares  containing $1's$ in $\beta(G)$ by $b_1, b_2, \ldots, b_{k-2}, f_1, f_2, \ldots, f_{t-q_k}, b_k$,   from left to right.
We also label the square $(p_{k-1}, q_{k-1})$ by $b_{k-1}$ (or $a_{k-1}$), see Figure \ref{case4} for example.

\begin{figure}[h,t]
\begin{center}
\begin{picture}(50,30)
\setlength{\unitlength}{7mm} \linethickness{0.4pt}
\put(0,0){\line(1,0){2}}\put(1,0){\circle*{0.15}}
\put(2,0){\line(0,1){1}}\put(2,1){\circle*{0.15}}\put(2,1){\line(1,0){1}}
\put(3,1){\line(0,1){1}}\put(3,2){\circle*{0.15}}\put(3,2){\line(1,0){1}}
\put(4,2){\line(0,1){1}}\put(4,3){\circle*{0.15}}\put(4,3){\line(1,0){3}}
\put(0,1){\line(1,0){1}}\put(1,1){\circle{0.15}}
\put(1,1){\line(0,1){1}}\put(1,2){\line(1,0){1}}\put(2,2){\circle{0.15}}
\put(2,2){\line(0,1){1}}\put(2,3){\line(1,0){1}}\put(3,3){\circle{0.15}}
\put(3,3){\line(0,1){1}}\put(3,4){\line(1,0){1}}\put(4,4){\circle{0.15}}
\put(4,4){\line(0,1){1}}\put(4,5){\line(1,0){5}}\put(6,5){\circle*{0.15}}\put(7,3){\line(0,1){3}}
\put(7,4){\circle*{0.15}}\put(8,3.7){\circle*{0.15}}\put(7,3.7){\circle{0.15}}\put(5,4){$E_{2}$}\put(5,4){$E_{2}$}\put(7.3,5.3){$E_{2}$}
\put(9,3.3){\circle*{0.15}}\put(8,3.3){\circle{0.15}}\put(9,0){\circle{0.15}}
\put(1.15,-0.3){$a_{1}$} \put(2.15,0.7){$a_{2}$}\put(3.15,1.7){$a_{3}$}\put(4.15,2.7){$a_{4}$}\put(5.6,5.2){$a_{5}$}\put(7.1,3.9){$a_{6}$}
\put(0.6,1.2){$b_{1}$}\put(1.6,2.2){$b_{2}$}\put(2.6,3.2){$b_{3}$}\put(3.6,4.2){$b_{4}$}\put(9.15,0){$b_{6}$}\put(6,5.2){$(b_{5})$}
\put(8.1,3.6){$e_{1}$}\put(7.05,3.4){$f_{1}$}\put(9.1,3.2){$e_{2}$}\put(8,3){$f_{2}$}

\end{picture}
\end{center}
\caption{The labelling of squares in $\phi_4(\pi)$ for $k=6$.} \label{case4}
\end{figure}

In $\phi_1(\pi)$ or $\phi_2(\pi)$, we denote by $E_1$ the union of the following four  parts of the board: the board that is above $a_1$ but below $b_1$ and to the left of $a_1$, the board that is above $a_{k-1}$ but below $b_{k-1}$, to the left of $c_1$ and to the right of $a_{k-1}$,   the union of the rectangles with corners $a_i$ and $b_{i+1}$ for $i=1,2,\ldots, k-2$, and the board that is above $a_k$ and to the right of $c_1$, see Figures \ref{case1} and \ref{case2} for example.

 We claim
that   there are no $1's$ inside  $E_1$ in  $\pi$,  $\phi_1(\pi)$ or
$\phi_2(\pi)$. By the selection  of $a_{k}$, there is no $1$ to the
right of $c_1$ and above $a_k$ in $\pi$, $\phi_1(\pi)$ or
$\phi_2(\pi)$.  Suppose that there is a $1$ in the rectangle with
corners $a_i$ and $b_{i+1}$ for $i=1, 2, \ldots, k-2$, then that 1
combining with the $1's$ positioned at $a_1, a_2, \ldots, a_i,
a_{i+2}, \ldots, a_k$ would form an $H_k$ in $\pi$, which
contradicts the selection of $a_{i+1}$. If there is a $1$ above
$a_1$ but below $b_1$, then that $1$, combining with the $1's$
positioned at $a_2, a_3, \ldots, a_k$ would form a $H_k$ in $\pi$,
which contradicts the selection of $a_1$. If there is a $1$   above
$a_{k-1}$ but below $b_{k-1}$ and to the left of $c_1$, then that
$1$, combining with the $1's$ positioned at $a_2, a_3, \ldots, a_k$
would form a $H_k$ in $\pi$, which contradicts the selection of
$a_1$. Thus,   all the $1's$ are to the left of $E_1$ or  to the
right of  $E_1$  in  $\pi$,  $\phi_1(\pi)$ or $\phi_2(\pi)$.

In $\phi_3(\pi)$ and $\phi_4(\pi)$, we denote by $E_2$ the union of
the following four parts of the board: the board that is above $a_1$
but below $b_1$ and to the left of $a_1$, the board that is above
and to the right of  $a_{k-2}$ but below $a_{k-1}$,  and to the left
of $a_k$, the union of the rectangles with corners $a_i$ and
$b_{i+1}$ for $i=1,2,\ldots, k-2$, and the board that is above
$a_{k-1}$ and to the right of $a_{k}$,  see Figures \ref{case3} and
\ref{case4} for example.

We claim that   there are no $1's$ inside $E_2$ in $\pi$,
$\phi_3(\pi)$ and $\phi_4(\pi)$.  By similar arguments  in $E_1$,
one can easily verify that there are no $1's$ inside the board that
is above $a_1$ but below $b_1$ and to the left of $a_1$,     the
union of the rectangles with corners $a_i$ and $b_{i+1}$ for
$i=1,2,\ldots, k-2$, and the board that is above $a_{k-1}$ and to
the right of $a_{k}$. It remains to show that there are no $1's$
inside the board that is above $a_{k-2}$ but below $a_{k-1}$ and to
the left of $a_k$. According to the definition of $Q_k$,   all of
the $1's$ between $a_{k-1}$ and $a_k$ are above $a_{k-1}$. This
implies that there are no $1's$ inside the board that is below and
to the right of $a_{k-1}$, and to the left of $a_k$. Now suppose
that there is a $1$ inside the rectangle with corners $a_{k-2}$ and
$a_{k-1}$.   Suppose that this $1$ is at the square $(\pi_g,g)$. If
$(\pi_g, g)$ is below $a_k$, then the $1's$ positioned at the
squares $a_2, a_3, \ldots, a_{k-2}, (\pi_g,g), a_{k-1}, a_k$ would
form a $Q_k$ in $\pi$, which contradicts the selection of $a_1$.  If
$(\pi_g, g)$ is above $a_k$, then we have two cases. If there exists
a  $j$ such that $g\leq j<q_k$ and $j\in \mathcal{D}(\pi)$, then the
$1's$ positioned at   $a_1, a_2, \ldots, a_{k-2}, (\pi_g,g),
a_{k-1}$ would form an $H_k$ in $\pi$, which contradicts the
selection of $a_{k-1}$. Otherwise,
  the $1's$ positioned at the
squares $a_1, a_2, \ldots, a_{k-2}, (\pi_g,g), a_k$ would form a
$Q_k$ in $\pi$, which contradicts the selection of $a_{k-1}$. Hence,
we have concluded that there are no $1's$ inside the board that is
above $a_{k-2}$ but below $a_{k-1}$ and to the left of $a_k$. Hence,
the claim is proved. In other words, all the $1's$ of $\pi$,
$\phi_3(\pi)$ or $\phi_4(\pi)$ are to the left or to the right of
$E_2$.

\begin{definition}
 A $1$ is said to be  strictly to the left (resp. right) of $E_1$(or
 $E_2$) if it is lying to the left (resp.  right) of $E_1$ (or $E_2$)
 and does not belong to the boundary of $E_1$ (or $E_2$).
\end{definition}

We need the following useful properties that will play an essential
role in the construction of  {\em vertical slide algorithm} and {\em
horizontal slide algorithm} for $\phi$
 \noindent{\bf Properties}
\begin{itemize}
\item[(1)] For any $1\leq i \leq k-2$, the board that is above $a_1$   and below $b_i$ cannot contain a $J_i$ with all its $1's$ strictly to the left of $E_1$ (or $E_2$) in   $\phi(\pi)$.

    \item[(2)] For any $1\leq i<j\leq  k-2$, the rectangle with corners $b_i$  and  $b_j$   cannot contain a $J_{j-i}$ with all its $1's$ strictly to the left of $E_1$ (or $E_2$) in $\phi(\pi)$.
        Moreover, the rectangle with corners $b_i$  and  $b_{k-1}$   cannot contain a $J_{k-1-i}$ with all its $1's$ strictly to the left of $E_1$ in   $\phi_1(\pi)$ (or $\phi_2(\pi)$).

\end{itemize}
\pf
\begin{itemize}
\item[(1)] If there is such a $J_i$ below $b_i$ in $\phi(\pi)$, then it is below $a_{i+1}$. Therefore these $i$ $1's$, combining with   $a_{i+1}, a_{i+2}, \ldots, a_k$, will either form an $H_k$ or a $Q_k$ in $\pi$, which contradicts the selection of $a_1$.

\item[(2)] If there is a  $ J_{j-i}$  in this region,
then either its leftmost $1$ is to the left to $b_{i+1}$ (and hence
to the left of $a_{i+1}$), or else it lies to the right of $b_{i+1}$
(and $a_{i+1}$). In the first case,    $a_1, a_2, \ldots, a_i$,
combining with this $J_{j-i}$ and    $a_{j+1}, \ldots, a_k$,   will
form an  $H_{k}$ (or $Q_k$) in $\pi$, which contradicts the
selection of $a_{i+1}$. In the second case,     $a_2, a_3, \ldots,
a_{i+1}$, combining with this $J_{j-i}$  and    $a_{j+1}, \ldots,
a_k$, will form an $H_{k}$ (or $Q_k$) in $\pi$, which contradicts
the selection of $a_{1}$. \qed
\end{itemize}

Now we proceed to introduce  the   vertical slide algorithm  and
 horizontal slide algorithm  for $\phi$.

Suppose that $H$ is   a  $J_t$ in $\phi(\pi)$. Label the squares containing  $1's$ of
$H$ by $h_1, h_2, \ldots, h_t$, from left to right.

 \noindent{\em  Vertical slide algorithm for $\phi$}: When $H$ is in $\phi_1(\pi)$ (or $\phi_2(\pi)$),  find the largest $i$ such that $b_i$ falls in $H$ with $i\leq k-1$; otherwise,  find the largest $i$ such that $b_i$ falls in $H$ with $i\leq k-2$.   If there is a
$1$ of $H$ which is below $b_i$ and to the right of $E_1$ (or
$E_2$), find the rightmost square containing  such a  $1$ and denote it by $h_y$.
 Find $x$
such that $h_y$ is to the right of $b_x$, and to the left of
$b_{x+1}$. By property $(2)$,  there are at most $i-x$ $1's$ in $H$,
which are above $b_x$ but not above $b_i$, and weakly  to the left
of $E_1$ (or $E_2$). So we can replace these $1's$   by    $b_{x+1},
b_{x+2}, \ldots, b_i$, and hence by   $a_{x+1}, a_{x+2}, \ldots, a_i$.

 We can repeat the vertical slide algorithm until the following two cases appear.\\
 (1) There is no $b_i$ that falls in $H$.   \\
 (2) There is such a $b_i$, but $h_y$ does not exist. By Property (1),
 there are at most $i$ $1's$  of $H$ that are above $a_1$ but not above $b_i$,  and  weakly to the left of $E_1$ (or $E_2$).  So
  we can replace these $1's$ by   $a_{1}, a_2, \ldots, a_i$ to form an  $J_t$ in $\pi$.

Suppose that $H$ is   a  $J_t$ in $\phi(\pi)$. Label the squares containing  the $1's$ of
$H$ by $h_1, h_2, \ldots, h_t$, from left to right. Assume that
$h_t$ is not above  $b_{k-1}$ when $H$ is in  $\phi_1(\pi)$ (or
$\phi_2(\pi)$), and not above $b_{k-2}$ when $H$ is in $\phi_3(\pi)$
(or $\phi_4(\pi)$).

 \noindent{\em  Horizontal  slide algorithm for $\phi$}:  When $H$ is in $\phi_1(\pi)$ (or $\phi_2(\pi)$),  find the least $i$ such that $b_i$ falls in $H$ with $i\leq k-1$; otherwise,  find the least $i$ such that $b_i$ falls in $H$ with $i\leq k-2$.
   If there is
a $1$ of $H$ which is above  $b_i$ and to the right of $E_1$ (or
$E_2$), find the leftmost square containing such a  $1$ and denote it by $h_y$. Find $x$ such
that $h_y$ is above  $b_x$, and below  $b_{x+1}$. By property $(2)$,
there are at most $x+1-i$ $1's$ in $H$, which are below  $b_{x+1}$
but not below  $b_i$, and weakly  to the left of $E_1$ (or $E_2$).
So we can replace these $1's$ by    $b_i, b_{i+1}, \ldots, b_x$, and
hence   by   $a_{i+1}, a_{i+2}, \ldots, a_{x+1}$.

 We can repeat the horizontal slide algorithm until the following two cases appear.\\
 (1) There is no $b_i$ that falls in $H$.  \\
 (2) There is such a $b_i$, but $h_y$ does not exist. Find the least
 $v$ such that $h_t$ is not above $b_v$.   By property $(2)$, we have at most $v-i$ $1's$ of
$H$ that are below $b_{v}$ but not below $b_i$ and  weakly to the
left of $E_1$(or $E_2$). So  we can replace these $1's$ by   $a_{i+1},
a_{i+2}, \ldots, a_{v}$ to form an $J_t$ in $\pi$.

Our next goal is to show  that the transformation $\phi$ have the
following properties, which are essential in the proof of Theorem
\ref{th3}.

\begin{lemma}\label{phid}
If there is no $F_k$ with at least one square in a row below $a_1$, then we have $\mathcal{D}(\pi)=\mathcal{D}(\phi(\pi))$.
\end{lemma}
 \pf Since there are no $1's$ inside $E_1$ (or $E_2$) and no $F_k$ with at least one square in a row below $a_1$, one can easily verify that $\mathcal{D}(\pi)=\mathcal{D}(\phi(\pi))$. The details are omitted here. \qed

\begin{lemma}\label{phi1}
If $\pi$ contains no $F_k$ with at least one square in a row below $a_1$, then $\phi(\pi)$ contains no such $F_k$.
\end{lemma}
\pf If not, suppose that $H$ is such an $F_k$ in $\phi(\pi)$. Label  the squares containing the $1's$ of  $H$ by $h_1, h_2, \ldots, h_k$, from left to right.
Then $h_k$ is below $a_1$.

 We claim that $h_k$ must be positioned  to the left of $a_{k-1}$.
 If not, then    $a_1, a_2, \ldots,$ $ a_{k-1},$ $ h_k$ would form an $F_k$ in $\pi$, which   contradicts the hypothesis.
 From the construction of the transformation  $\phi$, it follows that at least one of $b_1, b_2, \ldots, b_{k-2}$ must fall in $H$. Otherwise, $H$ is   an $F_k$ with   at least one square in a row below $a_1$ in $\pi$, which contradicts the  hypothesis.

By applying the vertical slide algorithm repeatedly to the $J_{k-1}$ consisting of    $h_1, h_2, \ldots, h_{k-1}$, we can get a $J_{k-1}$ not below $a_1$  in    $\pi$.  Then, that $J_{k-1}$ combining with     $h_k$ will form an $F_k$ in $\pi$, which  contradicts the hypothesis. \qed

\begin{lemma}\label{phi2}
The board that is to the left of $b_{t+1}$ and above $a_1$ cannot contain a $J_t$ in $\phi(\pi)$ with its highest $1$ below $b_t$ for $t=1,2,\ldots, k-1$.
\end{lemma}
\pf First  we aim to prove the assertion for $1\leq t\leq k-2$.   Suppose that $H$ is such a $J_t$ in $\phi(\pi)$.
 Label the squares containing the   $1's$  of $H$ by $h_1, h_2, \ldots, h_t$ from left to right. We claim that at least one of $ b_1, b_2, \ldots, b_{t-1}$ must
  fall in $H$. Otherwise, these $t$ $1's$, combining with  $a_{t+1}, a_{t+2}, \ldots, a_{k}$,  would form an $H_k$ or $Q_k$ in $\pi$.
  This  contradicts the selection of $a_1$.

By applying the   the horizontal slide algorithm repeatedly to $H$, we can get a
$J_{t}$
 in    $\pi$. It is easy to check that  the obtained $J_t$ is below and  to the left of $a_{t+1}$   and above $a_1$.
 That $J_t$, combining with   $a_{t+1}, a_{t+2}, \ldots, a_{k}$, would form an $H_k$ or $Q_k$ in $\pi$. This  contradicts the selection of $a_1$.
   Thus, we have concluded that the assertion  holds for $1\leq t\leq k-2$.

Now we proceed to show that the assertion also holds for $t=k-1$.  Suppose that $G$ is a $J_{k-1}$ in $\phi(\pi)$, which
is to the left of $b_k$ and below $b_{k-1}$. We label the squares containing  the $1's$ of $G$ by $g_1, g_2, \ldots, g_{k-1} $, from left to right. We have three cases.

  Case 1. $G$ is in $\phi_1(\pi)$.   By  repeating  the horizontal slide algorithm, we can get a $J_{k-1}$ in $\pi$, which is to the left of $b_k$ and above $a_1$. Moreover, since $s-1\in \mathcal{D}(\pi)$ and $\mathcal{D}(\pi)=\mathcal{D}(\phi(\pi))$, we have $s-1\in \mathcal{D}(\phi(\pi))$. Recall that $b_k$ is at column $s$.  Thus,
   the obtained $J_{k-1}$ combining with   $a_k$ would form an $H_k$ in $\pi$. This contradicts the selection of $a_1$.

 Case 2. $G$ is in $\phi_2(\pi)$.  If $g_{k-1}\neq f_i$,  we can get a $J_{k-1}$ in $\pi$ by repeating the horizontal slide algorithm, which is below  and   to  the left of $a_k$  and above $a_1$. We label  the squares containing the $1's$ of this $J_{k-1}$ by $m_1, m_2, \ldots, m_{k-1}$, from left to right. If  $m_{k-1}$ is below $e_1$,  then   $m_2, m_3, \ldots, m_{k-1}$,
 combining with    $a_k, e_1$, would form a $Q_k$ in $\pi$, which contradicts the selection of $a_1$. If $m_{k-1}$ is above $e_1$, then it is positioned to the left of $c_1$ in $\pi$ since all the $1's$ positioned at columns $s, s+1, \ldots, q_k-1$  form a $J_{q_k-s}$, and   $(\pi_{q_k-1}, q_k-1)$ is below $e_1=(\pi_{q_k+1}, q_k+1)$.  Moreover, since $s-1\in \mathcal{D}(\pi)$ and $\mathcal{D}(\pi)=\mathcal{D}(\phi(\pi))$, we have $s-1\in \mathcal{D}(\phi(\pi))$.  Recall that $c_1 $ is at column $s$. Thus,
     $m_1, m_2, \ldots, m_{k-1}, a_k$ will form an $H_k$, which contradicts the selection of $a_1$.

  If $g_{k-1}=f_i$ for some $i$, then  we can get a $J_{k-1}$ in $\pi$ by repeating the horizontal slide algorithm and replacing  $f_j's$ by $e_j's$ whenever $f_j$ falls in $G$.  Notice that the rightmost $1$ of the obtained $J_{k-1}$ is $e_i$.  If $i=1$, then this $J_{k-1}$ combining with   $a_k$ would form a $Q_k$ in $\pi$.  For $i>1$, this $J_{k-1}$ combining with  $e_{i-1}$ would form a $Q_k$ in $\pi$. In both cases, we get a  $Q_k$ that is above $a_1$. This contradicts the selection of $a_1$.

Case 3.  $G$ is in $\phi_3(\pi)$ or $\phi_4(\pi)$.  When
$g_{k-1}\neq f_i$, according to the definition of $Q_k$, there
contains no $1's$ which are  below and to the right of $a_{k-1}$,
and to the left of $a_k$.
  So $g_{k-1}$ is to the left of $a_{k-1}$.   Recall that we have shown that there is no $J_{k-2}$ in $\phi(\pi)$,
  which is to the left of $b_{k-1}$ and below $b_{k-2}$. So $g_{k-1}$ is   below and  to the left of $b_{k-1}$ (and $a_{k-1}$), and
   above and to the left of  $b_{k-2}$.
   By repeating the vertical slide algorithm, we can get a $J_{k-1}$ not below $a_1$ in $\pi$, whose rightmost $1$ is positioned at
   $g_{k-1}$. Then this $J_{k-1}$ combining with  $a_{k-1}$ would
   form an $H_k$
    in $\pi$,  which contradicts the selection of $a_{k-1}$.   When $g_{k-1}=f_i$ for some $i\geq 1$,   by the same way as in Case 2,    we can get a  $Q_k$  above $a_1$ in $\pi$. This contradicts the selection of $a_1$.

 Thus, we deduce  that  the assertion also holds for $t=k-1$. This completes the proof.\qed

\begin{lemma}\label{phi3}
The rows above $a_1$ cannot contain an $H_k$ or $Q_k$ in $\phi(\pi)$.
\end{lemma}

In order to prove Lemma \ref{phi3}, we need the following two lemmas.

\begin{lemma}\label{phi4}
Suppose that  $G$ is   an $H_k$  above  $a_1$  in $\phi(\pi)$. Label the squares containing the $1's$  of  $G$ by $g_1, g_2, \ldots, g_k$, from left to right. Then the squares $g_{k}$ and $g_{k-1}$ are also filled with $1's$ in $\pi$.
\end{lemma}
\pf    Here we only prove the assertion for  $\phi_1(\pi)$ and $\phi_4(\pi)$.  All the other cases can be verified by similar arguments. By Lemma \ref{phi2}, there is no $J_{k-1}$ below $b_{k-1}$ and to the left of $b_k$ in $\phi_1(\pi)$  and $\phi_4(\pi)$. This implies that neither $g_{k-1}$ nor  $g_k$  will be any of $b_i's$  for  $1\leq i\leq k-2$ in $\phi_1(\pi)$ and  $\phi_4(\pi)$.  Moreover, neither $g_{k-1}$ nor  $g_k$ will be any of $f_i's$ in $\phi_4(\pi)$.   Thus, we have deduced that the assertion holds for $\phi_4(\pi)$.

    In order to prove the assertion  for $\phi_1(\pi)$, it remains to show that neither $g_k$ nor $g_{k-1}$ will be any of $b_{k-1}$ and $d_i's$ in $\phi_1(\pi)$.  We have four cases.
    \begin{itemize}
    \item[(1)]If  $g_k=b_{k-1}$, then    $g_1, g_2, \ldots, g_{k-1}$ form   a $J_{k-1}$, which is  to the left of $b_{k}$ and below $b_{k-1}$ in $\phi_1(\pi)$. This  contradicts Lemma \ref{phi2}.
    \item[(2)]If  $g_{k-1}=b_{k-1}$, then $g_k$ is above $b_{k-1}$ and to the left of $E_1$. This implies  the square  $g_k$ is also filled with a $1$ in $\pi$.  Since $\mathcal{D}(\pi)=\mathcal{D}(\phi(\pi))$, the $1's$ positioned at $b_{k-1}$ and $g_k$ belong to two different blocks of $\phi(\pi)$ imply that those positioned at $a_{k-1}$ and $g_k$ also belong to two different blocks of $\pi$. Thus,
   $a_1, a_2, \ldots, a_{k-1}, g_{k}$   form an $H_k$ in $\pi$. This contradicts the selection of $a_k$ since $g_k$ is above $a_k$.
   \item[(3)]      If $g_k=d_i$ for some $i$,   then  we have that $g_{k-1}$ is to
   the left of $b_k$ since  $d_1$, $d_2$, $\ldots, d_{q_k-s}$ lie in consecutive columns and  form a $J_{q_k-s}$. Thus,   $g_1, g_2, \ldots, g_{k-1}$ will form a $J_{k-1}$ in $\phi_1(\pi)$, which is to the left of
   $b_k$ and below $b_{k-1}$. This contradicts Lemma \ref{phi2}. Hence, we have $g_k\neq d_i$ for any $i\geq 1$.
  \item[(4)] If $g_{k-1}=d_{j}$ for some $j\geq 1$, then    $g_k$ is to the right of $a_k$ since $d_1$, $d_2$, $\ldots, d_{q_k-s}$
  lie in consecutive columns and form a $J_{q_k-s}$.
   By repeating the horizontal slide algorithm for $\phi$ and replacing any of $d_i's$ that falls in $G$  by $c_i's$,
   we will get a $J_{k-1}$ in $\pi$ whose rightmost $1$ is $c_j$.  since $\mathcal{D}(\pi)=\mathcal{D}(\phi(\pi))$  and $q_k\in \mathcal{D}(\pi)$, we have $q_k\in \mathcal{D}(\phi(\pi))$. Thus, this
        $J_{k-1}$ combining with $g_k$ will form an $H_k$
     in $\pi$.  This contradicts the selection of $a_1$.
\end{itemize}
Hence,  we have concluded that the assertion  also holds for $\phi_1$, which completes the proof.\qed

\begin{lemma}\label{phi5}
Suppose that  $H$ is   a  $Q_k$  above  $a_1$  in $\phi(\pi)$, in
which the last two $1's$ lie in two consecutive columns. Label the squares containing  the
$1's$  of  $H$ by $h_1, h_2, \ldots, h_k$, from left to right. Then
 the squares $h_{k}$ and $h_{k-1}$ are also filled with $1's$ in
 $\pi$.
\end{lemma}
\pf
   Here we only prove the assertion for  $\phi_1$ and $\phi_4$.  All the other cases  can be verified by similar arguments. By Lemma \ref{phi2}, there is no $J_{k-1}$ below $b_{k-1}$ and to the left of $b_k$ in $\phi_1(\pi)$.   This implies that neither $h_{k-1}$ nor  $h_k$  will be any of $b_i's$  for  $1\leq i\leq k-2$ in $\phi_1(\pi)$ and  $\phi_4(\pi)$.  Moreover, neither $g_{k-1}$ nor  $g_k$ will be any of $f_i's$ in $\phi_4(\pi)$.   Thus, we deduce that the assertion holds for  $\phi_4$.

In order to prove the assertion  for $\phi_1$, it remains to show that neither $h_k$ nor $h_{k-1}$ will be any of $b_{k-1}$ and $d_i's$ in $\phi_1(\pi)$.
 Since $b_k, d_1, d_2, \ldots, d_{q_k-s}$  lie in consecutive columns and  form a $J_{q_k-s+1}$  in $\phi_1(\pi)$, neither
 of $d_i's$ can be $h_k$. Moreover, neither of $d_i's$ can be $h_{k-1}$ for $1\leq i\leq q_k-s-1$. Thus we have   $h_{k-1}=d_{q_k-s}$,  $  h_{k-1}=b_{k-1}$ or $h_k=b_{k-1}$.

   If $h_{k-1}=d_{q_k-s}$, then  by applying  the horizontal slide algorithm  for $\phi$ repeatedly to   $h_1, h_2, \ldots, h_{k-2}, h_k$ and replacing any of $d_i's$ that falls in $h_1, h_2, \ldots, h_{k-2}, h_k$  by $c_i's$, we will get a $J_{k-1}$ above $a_1$ in $\pi$. Notice that the rightmost $1$ of the obtained $J_{k-1}$ is $h_k$.  This $J_{k-1}$, combining with   $a_k$,
   will form a $Q_k$ in $\pi$, which contradicts the selection of $a_1$.
  If $ h_k=b_{k-1}$,  then    $a_1, a_2, \ldots, a_{k-2}$, combining with  $h_{k-1}$ and $a_{k-1}$, will form a $Q_k$ in $\pi$, which contradicts the selection of $a_{k-1}$.
   If  $ h_{k-1}=b_{k-1}$, then $h_k$ is below  $a_{k-1}$ and to the left of $b_k$. Then    $h_1, h_2, \ldots,h_{k-2}, h_{k}$ form a $J_{k-1}$ in $\phi_1(\pi)$, which is to the left of $b_k$ and below $b_{k-1}$. This contradicts Lemma \ref{phi2}.
           Hence, we have proved that the assertion  also holds for $\phi_1$.  \qed

\noindent{\bf The proof of Lemma \ref{phi3}.}  If not, suppose that    $G$ is   an $H_k$  above  $a_1$  in $\phi(\pi)$. Label the squares containing  $1's$ of  $G$ by $g_1, g_2, \ldots, g_k$, from left to right.   Moreover, let $H$ be a $Q_k$ above $a_1$ in $\phi(\pi)$ such that the rightmost two $1's$ lie in two consecutive columns.  We label the squares containing  the  $1's$ of  $H$ by $h_1, h_2, \ldots, h_k$.  According to the definition of $Q_k$,    there is a  $Q_k$ above $a_1$ in $\phi(\pi)$ if and only if  there  exists such an $H$.

 We wish to replace some $1's$ of $G$ (resp.  $H$ ) to form an $H_k$ (resp. $Q_k$) in $\pi$. Here we only consider the case
 when $G$ (resp. $H$) is in $\phi_1(\pi)$.
 The other cases can be verified by the similar arguments.  Since the transformation $\phi_1$  does not change the positions of any other $1's$,
  one of $b_i's$ and $d_i's$  must fall in $G$ (resp.  $H$).

   First, replace each $d_i$ by $c_i$ whenever $d_i$ falls in $g_1, g_2, \ldots, g_{k-2}$ (resp. $h_1, h_2, \ldots, h_{k-2}$).
 Then, find the largest $i$ such that  $b_i$ falls in $g_1, g_2, \ldots, g_{k-2}$ (resp. $h_1, h_2, \ldots, h_{k-2}$).
 We can apply the vertical slide algorithm repeatedly to   $g_1, g_2, \ldots, g_{k-2}$ (resp. $h_1, h_2, \ldots, h_{k-2}$) until the following two cases appear.\\
 (1) There is no $b_i$ that falls in $g_1, g_2, \ldots, g_{k-2}$ (resp. $h_1, h_2, \ldots, h_{k-2}$) . \\
 (2) There is such a $b_i$, but there is no $1$ positioned at  the squares  $g_1, g_2, \ldots, g_{k-2}$ (resp. $h_1, h_2, \ldots, h_{k-2}$) that is to the left of $b_i$ and to the right of $E_1$. Since there are at most $i$ $1's$   positioned at   $g_1, g_2, \ldots, g_{k-2}$ (resp. $h_1, h_2, \ldots, h_{k-2}$)  that are not above $b_i$ and    to the left of $E_1$, we can replace these $1's$ by $a_{1}, a_2, \ldots, a_i$.

 In both cases, we get a $J_{k-2}$    not below $a_1$ in $\pi$. From Lemmas \ref{phi4} and $\ref{phi5}$,  the squares $g_k$ and $g_{k-1}$ (resp. $h_k$ and
 $h_{k-1}$) are also filled with $1's$ in $\pi$.
 Recall that   $\mathcal{D}(\pi)=\mathcal{D}(\phi(\pi))$ and the $1's$ positioned at $g_{k-1}$ and $g_k$ belong to two different blocks of $\phi(\pi)$. This yields that   the $1's$ positioned at $g_{k-1}$ and $g_k$ also belong to two different blocks of $\pi$.  Thus,    the obtained $J_{k-2}$, combining with  $g_{k-1}$ and $g_k$ (resp. $h_{k-1}$ and $h_k$)  forms an $H_k$ (resp.  $Q_k$) in $\pi$.
 In the first case, the obtained $H_k$ (resp. $Q_k$) is above $a_1$,
 which contradicts the selection of $a_1$.
  In the second case, suppose that $g_z$(resp. $h_z$) is the first square containing a $1$ of the obtained $H_k$ (resp. $Q_k$) that is to the right of
  $a_i$. clearly,  $g_z$ (resp. $h_z$) is above $b_i$ and $a_{i+1}$.
    If $g_z$ (resp. $h_z$) is to the left of $a_{i+1}$, then the obtained $H_k$ (or $Q_k$)
      contradicts the selection of $a_{i+1}$. Otherwise,   $a_2, a_3, \ldots,  a_{i+1} $, combining with the $1's$ of the obtained $H_k$ (resp. $Q_k$) that are to the right of $a_i$, would form an $H_k$ (or $Q_k$ ) in $\pi$. This contradicts the selection of $a_1$, which  completes the proof. \qed

\subsection{The map $\Psi$  from the set $\mathcal{S}_n(H_k, Q_k)$ to the set $\mathcal{S}_n(F_k)$}

Before we describe  the map $\Psi$ we define three transformations, which will play an essential role in the construction of the map $\Psi$.

Let $\sigma=\{(\sigma_1, 1), (\sigma_2, 2), \ldots, (\sigma_n,n)\}$.
 Suppose that  $G$ is the submatrix  of $\sigma$ at columns $c_1<c_2<\ldots<c_{k}<c_k+1<c_k+2<\ldots<t$ and rows $ r_1< r_2< \ldots<r_{k-1} > r_k< \sigma_{c_k+1}<\sigma_{c_k+2}  < \ldots< \sigma_t$, in which the squares $(r_i, c_i)$ are filled with $1's$ for all $i=1, 2, \ldots, k$. Let $\delta(G)$ be the submatrix at the same rows and columns as $G$,  such that the squares
    $(r_k, c_1)$, $(r_1, c_2)$, $\ldots$, $(r_{k-2}, c_{k-1})$, $(\sigma_{c_k+1}, c_k)$, $(\sigma_{c_k+2}, c_k+1)$, $\ldots$,  $ (\sigma_{t}, t-1)$, $(r_{k-1}, t)$ are filled with $1's$ and all the other squares are filled with $0's$.

Suppose that $H$ is the submatrix of $\sigma$ at  columns $c_1<c_2<\ldots<c_{k-1}< t<t+1<\ldots<c_k$ and rows $r_1< r_2<\ldots<r_{k-1}>\sigma_{t}>\sigma_{t+1}> \ldots> \sigma_{c_k}=r_k$,   in which the squares $(r_i, c_i)$ are filled with $1's$ for all $i=1, 2, \ldots, k$.   Define $\gamma(H)$ to be the submatrix at the same  columns  and  rows as $H$,  such that the squares
    $(r_k, c_1)$, $(r_1, c_2)$, $\ldots$, $(r_{k-2}, c_{k-1})$, $(r_{k-1}, t)$,     $(\sigma_{t}, t+1)$, $(\sigma_{t+1}, t+2)$, $\ldots$, $(\sigma_{c_k-1}, c_k)$
are filled with $1's$ and all the other squares are filled with $0's$.

{\noindent \em The transformation $\psi$:}
Suppose that $\sigma=\sigma_1\sigma_2\ldots \sigma_n $ is a permutation in $\mathcal{S}_n$.
  First, find the  lowest square  $(p_k, q_k)$ containing a $1$,  such that there is an $F_k$ in $\sigma$ in which the $1$ positioned at $(p_k, q_k)$ is its rightmost  $1$. Then,
find the lowest square  $(p_{k-1}, q_{k-1})$ containing a $1$,  such that there is an $F_k$ in $\sigma$ in which the $1'$ positioned at $(p_k, q_k)$  and $(p_{k-1}, q_{k-1})$ are the rightmost  two $1's$.
  Find $(p_{k-2}, q_{k-2}), (p_{k-3}, q_{k-3}), \ldots, (p_1, q_1)$ one by one as $(p_{k-1}, q_{k-1})$.  Assume  that there is no $H_k$ or $Q_k$ above row $p_k$ in $\sigma$.

   If $\sigma_{q_k-1}>\sigma_{q_k+1}$,   then  we wish to  generate a permutation $\pi$ from $\sigma$ by the considering the following two  cases.
 \begin{enumerate}
 \item[Case 1.]   $p_{k-1}=\sigma_{q_{k-1}}<\sigma_{q_{k-1}+1}<\ldots<\sigma_{q_k-1}>\sigma_{q_k}=p_k$  and
  $ \sigma_{q_k+1}>p_{k-2}$. In this case, let $G$ be the submatrix of $\sigma$ at columns $q_1<q_2<\ldots<q_{k-2}<q_k$ and rows $p_k<p_1<p_2<\ldots<p_{k-2}$.    Replace $G$ by $\theta^{-1}(G)$ and leave all the other rows and columns fixed.
 \item[Case 2.] Otherwise,  find the least $t$    such that $t>q_k$ and  $t\in \mathcal{D}(\sigma)$. If such $t$ does not exist, set $t=n$. In this case, let  $G$ be the submatrix of $\sigma$ at columns $q_1<q_2<\ldots<q_k<q_k+1<\ldots<t$ and rows $p_1<p_2<\ldots p_{k-1}>p_k<\sigma_{q_k+1}<\sigma_{q_k+2}<\ldots<\sigma_t$.    Replace $G$ by $\delta(G)$  the other rows and columns fixed.
 \end{enumerate}
   If $q_k=n$  or $\sigma_{q_k-1}<\sigma_{q_k+1}$,   then we wish to  generate a permutation $\pi$ from $\sigma$ by   considering the following two  cases.
 \begin{enumerate}
 \item[Case 3.] If  there exists an $s$ such that $q_{k-1}<s<q_k$ and $\sigma_{s-1}>\sigma_s<\sigma_{s+1}$. Find the largest $t$ such that $q_{k-1}<t\leq q_k$ and $t-1\in \mathcal{A}(\sigma)$. Let $G$ be the submatrix of $\sigma$ at columns $q_1<q_2<\ldots<q_{k-1}<t<t+1<\ldots<q_k$ and rows $p_1<p_2< \ldots<p_{k-1}>\sigma_{t}>\sigma_{t+1}>\ldots>p_k $. Replace $G$ with $\gamma(G)$ and leave all the other rows and columns fixed.

 \item[Case 4.] Otherwise, we have $p_{k-1}=\sigma_{q_{k-1}}<\sigma_{q_{k-1}+1}<\ldots<\sigma_{t-1}>\sigma_t>\sigma_{t+1}>\ldots >\sigma_{q_k}=p_k$  for some   $t$  with $q_{k-1}<t\leq q_k$. Let $G$ be the submatrix of $\sigma$ at columns $q_1<q_2<\ldots<q_{k-2}<q_{k-1}<t<t+1<\ldots<q_k$ and rows $p_1<p_2< \ldots<p_{k-2}<p_{k-1}>\sigma_{t}>\sigma_{t+1}>\ldots>p_k $. Replace $G$ with $\gamma(G)$ and leave all the other rows and columns fixed.

  \end{enumerate}

\begin{remark}
In Case $2$, the selection of $(p_k,q_k)$ ensures that $p_k<\sigma_{q_k+1}$. If not,  the the $1's$ positioned at $(p_1, q_1), (p_2, q_2), \ldots, (p_{k-1}, q_{k-1}), (\sigma_{q_k+1}, q_k+1)$ would form an $F_k$, which contradicts the selection of $(p_k,q_k)$.
In Case $3$, the existence of such $s$ and the hypothesis that  there is no $H_k$ above row $p_k$ ensure  that $p_{k-1}>\sigma_t$. If not, then the $1's$ positioned at $(p_2, q_2), (p_3, q_3), \ldots, (p_{k-1}, q_{k-1}), (\sigma_t, t)$ would form an $H_k$ above row $p_k$ in $\sigma$.  In Case $4$,   the hypothesis that  there is no $Q_k$ above row $p_k$ ensures  that $p_{k-1}>\sigma_t$.  If not, then the $1's$ positioned at $(p_2, q_2), (p_3, q_3), \ldots, (p_{k-1}, q_{k-1}), (\sigma_{t-1}, t-1), (\sigma_{t}, t)$ would form a  $Q_k$ above row $p_k$ in $\sigma$.
\end{remark}

\begin{remark}
We denote the resulting permutation in   Case $1$,    Case $2$, Case $3$ and Case $4$   by $\psi_1(\sigma)$, $\psi_2(\sigma)$, $\psi_3(\sigma)$ and $\psi_4(\sigma)$,  respectively.
\end{remark}

It is obvious that the transformation $\psi$   changes every occurrence of $F_k$ to an occurrence of $H_k$ (or  $Q_k$).
 Denote by    $\Psi$   the iterated transformation, that recursively transforms every occurrence of $F_k$ into  $H_k$ (or  $Q_k$).

Using the notation of the algorithm for $\psi_1$, we label the
squares   containing $1's$ in $G$ by $b_1, b_2, \ldots, b_{k}$, and
the squares containing $1's$ in $\theta^{-1}(G)$ by $a_1, a_2,
\ldots, a_{k-1}, a_k$, from left to right, see Figure \ref{case3}
for example.

Using the notation of the algorithm for $\psi_2$, we label the squares containing $1's$ in $G$ by $b_1, b_2, \ldots, b_{k-1},  b_k, d_1, d_2, \ldots, d_{t-q_k} $, and the squares containing $1's$ in $\delta(G)$ by $a_1, a_2, \ldots, a_{k-1},  c_1, c_2, \ldots, c_{t-q_k}, a_k$, from left to right, see Figure \ref{case1} for example.

Using the notation of the algorithm for $\psi_3$, we label the
squares containing $1's$ in $G$ by $b_1, b_2, \ldots, b_{k-1},  f_1,
f_2, \ldots, f_{q_k-t}, b_{k} $, and the squares containing $1's$ in
$\gamma(G)$ by $a_1, a_2, \ldots, a_{k-1}, a_k, e_1, e_2, \ldots,
e_{q_k-t}$, from left to right. We also label the minimum entry of
the block to which $f_1$ belongs  by  $c_1$,  see Figure \ref{case2}
for example.

Using the notation of the algorithm for $\psi_4$, we label the squares   containing $1's$ in $G$ by $b_1, b_2, \ldots, b_{k-1},  f_1, f_2, \ldots, f_{q_k-t}, b_{k} $, and the squares containing $1's$ in $\gamma(G)$ by $a_1, a_2, \ldots, a_{k-1}, a_k, e_1, e_2, \ldots, e_{q_k-t}$, from left to right, see Figure \ref{case4} for example.

In $\psi_2(\sigma)$ and  $\psi_3(\sigma)$, let $E_1$ be the same board defined in $\phi_1(\pi)$ and $\phi_2(\pi)$. Similarly, in $\psi_1(\sigma)$ and $\psi_4(\sigma)$,  let $E_2$ be the same board defined in $\phi_3(\pi)$ and $\phi_4(\pi)$. From  the selection of $b_i's$ and the hypothesis that there is no $H_k$ above $a_1$,  it follows that there are no $1's$ inside $E_1$ (or $E_2$). In other words, all the $1's$ are to the left  or to the right of $E_1$ (or $E_2$) in $\psi(\pi)$.

Now we proceed to prove that the transformation $\psi$ have the
following properties, which are essential in the proof of Theorem
\ref{th3}.

 \begin{lemma}\label{psi0}
If there is no $H_k$ or $Q_k$ above $a_1$ in $\sigma$, then we have  $\mathcal{D}(\sigma)=\mathcal{D}(\psi(\sigma))$.
\end{lemma}
\pf  Since there are no $1's$ inside $E_1$ (or $E_2$) and no $H_k$  or $Q_k$ above  $a_1$, one can easily verify that $\mathcal{D}(\pi)=\mathcal{D}(\psi(\sigma))$. The details are omitted here. \qed

\noindent{\bf Properties}
 \begin{itemize}
    \item[($1'$)] For any $1\leq i<j\leq  k-1$, the rectangle with corners $a_i$  and  $a_j$   cannot contain a $J_{j-i}$ with all its $1's$
     strictly to the right of $E_1$ (or $E_2$) in $\psi(\sigma)$.
     \item[($2'$)]
    For any $1\leq i\leq  k-2$,  the rectangle with corners $a_i$  and  $a_{k}$   cannot contain a $J_{k-i-1}$ with all its $1's$
     strictly to the right of   $E_2$  in $\psi_1(\sigma)$ (or $\psi_4(\sigma)$).

\end{itemize}
\pf \begin{itemize}
\item[$(1')$]  If there is a  $ J_{j-i}$  in this region, then     $b_1, b_2,
\ldots, b_{i-1}$, combining with this $J_{j-i}$ and   $ b_j,  b_{j+1},
b_{j+2},  \ldots, b_k$,  will form an  $F_k$  in $\sigma$, which
contradicts the selection of $b_{j-1}$.

\item[$(2')$]  If there is a  $ J_{k-1-i}$  in this region, then the rightmost $1$ of this $J_{k-1-i}$
is to the left of $b_{k-1}$ since all the $1's$ lying between $b_{k-1}$ and $a_k$ are to the left of $E_2$.
 Clearly, the rightmost $1$ of this $J_{k-1-i}$ is below $b_{k-2}$.   So     $b_1, b_2,
\ldots, b_{i-1}$, combining with this $J_{k-1-i}$ and    $b_{k-1},
b_{k}$, will form an $F_k$   in $\sigma$, which contradicts the
selection of $b_{k-2}$.
\end{itemize}

\begin{lemma}\label{psi1}
  $\psi(\sigma)$ contains no $F_k$ with at least one square in a row below $a_1$.
\end{lemma}
\pf If not, suppose $H$ is such an $F_k$ in $\psi(\sigma)$. Label the squares containing the $1's$ of $H$ by $h_1, h_2,
\ldots, h_k$ from left to right. Then $h_k$ is below $a_1$. As in the proof of Lemma \ref{phi1}, we shall replace some $1's$ of $H$ (except $h_k$)
to form an $F_k$ in $\pi$, which contradicts the selection of $b_k$.

By the selection of $b_k$, we have that $h_k$ must be at the
left  side of $b_{k-1}$. From the construction of $\psi(\sigma)$, at
least one of $a_1, a_2, \ldots, a_{k-2}$ must fall in $H$.
Otherwise, $H$ is also an $F_k$ in $\sigma$, which contradicts the
selection of $b_k$.

Find the least $i$ such that $a_i$ falls in $H$.

\noindent{\em Vertical slide algorithm for $\psi$:} If there is a $1$ of $H$ which is above $a_i$ and to the left of $E$, find the leftmost square containing such a   $1$
 and denote  it by $h_y$. find $x$ such that $h_y$ is to the right of $a_x$ and to the left of $a_{x+1}$. Then by property
 ($1'$), there are at most $x+1-i$ $1's$ in $H$ that are below
 $a_{x+1}$ but
not below $a_i$, and to the right of $E_1$ (or $E_2$). So we can
replace these $1's$ by   $a_i, a_{i+1}, \ldots, a_x$, and hence by those positioned at
$b_i, b_{i+1}, \ldots, b_x$.

We can repeat the vertical slide algorithm until one of the following two cases appears.
\begin{itemize}
\item[(1)] There is no $a_i$ that falls in $H$. This ends the proof.
\item[(2)] There is such an $a_i$, but $h_y$ does not exist. Then suppose $a_v$ is the first square
 to the right of $h_k$. By property ($1'$) , there are at most $v-i$ $1's$ in $H$ that are below
 and to the left of $a_v$, but not below $a_i$,  and to the right of $E_1$ (or $E_2$).
 So we can replace these $1's$ by  $a_i, a_{i+1}, \ldots, a_{v-1}$, and hence by  $b_i, b_{i+1}, \ldots, b_{v-1}$. Then we have an $F_k$ in $\sigma$ with a square $h_k$ below $a_1$. \qed
\end{itemize}

\begin{lemma}\label{psi2}
If $\sigma$ contains no $H_k$ or $Q_k$ that is above $a_1$, neither does $\psi(\sigma)$.
\end{lemma}

In order to prove Lemma \ref{psi2}, we need the following two lemmas.

\begin{lemma}\label{psi3}
Suppose that  $G$ is   an $H_k$  above  $a_1$  in $\psi(\sigma)$. Label the squares containing the $1's$ of  $G$ by $g_1, g_2, \ldots, g_k$, from left to right. If $\sigma$ contains no $H_k$ or $Q_k$ that is above $a_1$,  then the squares $g_{k}$ and $g_{k-1}$ are also filled with $1's$ in $\sigma$.
\end{lemma}

\begin{lemma}\label{psi4}
Suppose that  $H$ is   a  $Q_k$  above  $a_1$  in $\psi(\sigma)$, in
which the last two $1's$ lie in two consecutive columns. Label the squares containing the
$1's$  of  $H$ by $h_1, h_2, \ldots, h_k$, from left to right.
 If
$\sigma$ contains no $H_k$ or $Q_k$ that is above $a_1$,  then the
squares $h_{k}$ and $h_{k-1}$ are also filled with $1's$ in
$\sigma$.
\end{lemma}

 Before we prove Lemmas \ref{psi3} and \ref{psi4}, we introduce the following {\em horizontal slide algorithm }  for $\psi$.

Suppose $H$ is   a  $J_k$ in $\psi(\sigma)$. Label the squares containing the $1's$ of  $H$ by $h_1, h_2,
\ldots, h_k$ from left to right.

\noindent{\em Horizontal slide algorithm for $\psi_2$ (or
$\psi_3$):} Find the largest $i$ such that $a_i$ falls in $H$ with
$i\leq k-1$.  If there is a $1$ of $H$ which is below $a_i$ to the
left of $E_1$, find the rightmost squares containing such a  $1$ and denote it by $h_y$. Find
$x$ such that $h_y$ is below $a_x$, and above $a_{x-1}$. Then by
property $(1')$, there are $i-x+1$ $1's$ in $H$ that are above
$a_{x-1}$ but not above $a_{i}$, and to the right of $E_1$. So we
can replace these $1's$ by   $a_{x}, a_{x+1}, \ldots, a_i$, and hence
by  $b_{x-1}, b_x, \ldots, b_{i-1}$.

We can repeat this horizontal slide algorithm until one of the following two cases appears.
\begin{enumerate}
\item[(1)] There is no $a_i$ that falls in $H$.
\item[(2)] There is such an $a_i$, but $h_y$ does not exist.
Find $x$ such that $h_1$ is below $a_{x+1}$ and above $a_x$. Then by property $(1')$, there are $i-x$ $1's$ in $H$ that are above $a_{x}$ but
not above $a_{i}$.  So we can replace these $1's$ by   $a_{x+1}, a_{x+2}, \ldots, a_i$, and hence by $b_{x}, b_{x+1}, \ldots, b_{i-1}$.
\end{enumerate}

\noindent{\em Horizontal slide algorithm for $\psi_1$ (or
$\psi_4$).} Find the largest $i$ such that $a_i$ falls in $H$ with
$i\leq k-2$ or $i=k$.  If there is a $1$ of $H$ which is below $a_i$
to the left of $E_2$, find the rightmost square containing such a $1$ and denote it by $h_y$.
Find $x$ such that $h_y$ is below $a_x$, and above $a_{x-1}$.  If
$i\leq k-2$, then by property $(1')$, there are $i-x+1$ $1's$ in $H$
that are above $a_{x-1}$ but not above $a_{i}$, and to the right of
$E_2$. So we can replace these $1's$ by those positioned at $a_{x}, a_{x+1}, \ldots,
a_i$, and hence by $b_{x-1}, b_x, \ldots, b_{i-1}$. If $i= k$, then
by property $(2')$, there are $k-x$ $1's$ in $H$ that are above
$a_{x-1}$ but not above $a_{k}$, and to the right of $E_2$. So we
can replace these $1's$ by    $a_{x}, a_{x+1}, \ldots, a_{k-2},
a_{k}$, and hence by $b_{x-1}, b_x, \ldots, b_{k-2}$.

We can repeat this horizontal slide algorithm until one of the following two cases appears.
\begin{enumerate}
\item[(1)] There is no $a_i$ that falls in $H$.
\item[(2)] There is such an $a_i$, but $h_y$ does not exist.
 Find $x$ such that $h_1$ is below $a_{x+1}$ and above $a_x$.
 If $i<k-2$, then by property $(1')$, there are $i-x$ $1's$ in $H$ that are above $a_{x}$ but  not above $a_{i}$.
   So we can replace these $1's$ by   $a_{x+1}, a_{x+2}, \ldots, a_i$, and hence by $b_{x}, b_{x+1}, \ldots, b_{i-1}$.
   If $i=k$, then by property $(2')$, there are $k-1-x$ $1's$ in $H$ that are above $a_{x}$ but  not above $a_{k}$.
   So we can replace these $1's$ by    $a_{x+1}, a_{x+2}, \ldots, a_{k-2},  a_k$, and hence by $b_{x}, b_{x+1}, \ldots, b_{k-2}$.
\end{enumerate}

\noindent{\bf The proof of Lemma \ref{psi3}.}  Here we only prove
the assertion for  $\psi_2(\sigma)$ and $\psi_4(\sigma)$. The other
cases can be verified by similar arguments. In order to prove the
assertion, it suffices to show that neither $g_k$ nor $g_{k-1}$ will
be any of the $a_i's$ and $c_i's$ in $\psi_{2}(\sigma)$, and be any
of  the $a_i's$ for $i=1,2,\ldots, k-2, k$ and $e_i's$ in
$\psi_4(\sigma)$.

We claim  there is no $J_{k-1}$ which is below $b_{k-1}$ but above
$a_1$, and  not to the right of $b_{k}$ in $\psi_2(\sigma)$ (or
$\psi_4(\sigma)$). If not, suppose that $R$ is such a $J_{k-1}$.
When $J_{k-1}$ is in $\psi_2(G)$, we can get a $J_{k-1}$ from $R$ by
repeating the horizontal slide algorithm for $\psi_2$.  When
$J_{k-1}$ is in $\psi_4(G)$, we can get a $J_{k-1}$ from $R$ by
repeating the horizontal slide algorithm for $\psi_4$ and replacing
any $e_i$ by $f_i$ whenever $e_i$ fall in $R$. In both cases,
  the obtained $J_{k-1}$ is below $b_{k-1}$ but above $a_1$,
and  to the left of $b_k$. Then Then this $J_{k-1}$ combining with
$b_k$ will form an $F_k$ in $\sigma$, which contradicts the
selection of $b_{k-1}$. Hence, the claim is proved.

From the  claim, it follows that     neither $g_k$ nor $g_{k-1}$
will be any of the $a_i's$  in $\psi_{2}(\sigma)$ for $ i\leq k-1$.
  In order to prove the assertion  for $\psi_2(\sigma)$, it remains to show that neither $g_k$ or $g_{k-1}$ will be any of $a_k$ and $c_i's$ in $\psi_2(\sigma)$.  Clearly, $g_{k-1}$ cannot be $a_k$ since there is no $1's$ above and to the right of $a_k$.
    \begin{itemize}
    \item[(1)]If  $g_k$ is either $a_k$ or one of $c_i's$ in $\psi_2(\sigma)$, then
    $g_1, g_2, \ldots, g_{k-1}$ form a $J_{k-1}$ which is to the left of $b_{k}$ and below $b_{k-1}$
     in $\phi_2(\sigma)$ since $c_1, c_2, \ldots, c_{t-q_k}, a_k$ lie in consecutive columns and form a
     $J_{t-q_k+1}$. This contradicts the claim proved above.

    \item[(2)]If  $g_{k-1}$ is one of $c_i's$ in $\psi_2(\sigma)$, then
    $g_k$ is to the right of $a_k$   since $c_1, c_2, \ldots, c_{t-q_k},$ $ a_k$ lie
    in consecutive columns and form a $J_{t-q_k+1}$. By repeating the horizontal
    slide algorithm for $\psi_2$ and replacing any $c_i$ falling in $G$ by $d_i$, we can get a $J_{k}$ above $a_1$ in $\sigma$ from $G$.
      Notice that if $g_{k-1}=c_j$, then the rightmost two $1's$ of the obtained $J_k$ are $g_k$ and $d_{j}$.  Recall that $a_k$   is positioned at column $t$.
      Since $\mathcal{D}(\sigma)=\mathcal{D}(\psi(\sigma))$  and $t\in \mathcal{D}(\sigma) $,  we have $t\in \mathcal{D}(\psi(\sigma))$. Thus,  the obtained $J_k$ is an $H_k$.
       This contradicts the fact that there contains no $H_k$ above $a_1$ in $\sigma$.

\end{itemize}
Hence,  we have concluded that the assertion    holds for
$\psi_2(\sigma)$.

From the claim that  is no $J_{k-1}$ which is below $b_{k-1}$ but
above $a_1$, and not to the right of $b_{k}$ in $\psi_4(\sigma)$,   it
follows that  neither $g_k$ nor $g_{k-1}$ will be any of the $a_i's$
for $i=1,2,\ldots, k-2, k$ and $e_i's$ in $\psi_{1}(\sigma)$. Hence,
we deduce that the assertion also holds for $\psi_4(\sigma)$, which
completes the proof.
   \qed

\noindent{\bf The proof of Lemma \ref{psi4}.}
   Here we only prove the assertion for  $\psi_2(\sigma)$ and $\psi_4(\sigma)$. The other cases can be verified by similar arguments.
In order to prove the assertion, it suffices to show that neither $h_k$ nor $h_{k-1}$ will be any of the $a_i's$ and $c_i's$ in $\psi_{2}(\sigma)$, and be any of  the $a_i's$ for $i=1,2,\ldots, k-2, k$ and $e_i's$ in $\psi_4(\sigma)$.

Recall that we have proved the claim in the proof of Lemma
\ref{psi3} that there is no $J_{k-1}$ which is below $b_{k-1}$ but
above $a_1$, and not to the right of $b_{k}$ in  $\psi_4(\sigma)$. It
follows that neither $h_k$ nor $h_{k-1}$ will be any of the $a_i's$
for $i=1,2,\ldots, k-2, k$ and $e_i's$ in $\psi_4(\sigma)$.
 Thus, the assertion holds for
$\psi_4(\sigma)$.

Similarly, from the claim proved in the proof of Lemma \ref{psi3},
it follows that neither $h_k$ nor $h_{k-1}$ will be any of the
$a_i's$ for $i=1,2,\ldots, k-1$ in $\psi_2(\sigma)$. In order to
prove the assertion for $\psi_2(\sigma)$, it remains to verify that
neither $h_k$ nor $h_{k-1}$ will be any of $a_k$ and $c_i's$ in
$\psi_2(\sigma)$. Recall that $c_1, c_2, \ldots, c_{q_k-t}, a_k$ lie
in consecutive columns and form a $J_{q_k-t+1}$ in $\psi_1(\sigma)$.
It implies that if $h_k$ or $h_{k-1}$ is one of $c_1, c_2, \ldots,
c_{q_k-t}, a_k$, then we have either $h_{k-1}=a_k$ or $h_{k}=c_1$.

In the former case,  we can get a $J_{k-1}$ above $a_1$ in $\sigma$
from the $J_{k-1}$  consisting of $h_1, h_2, \ldots, h_{k-2}, h_k$,
by repeating the horizontal slide algorithm for $\psi_2$ and
replacing any $c_i$   by $d_i$.  Since
$\mathcal{D}(\sigma)=\mathcal{D}(\psi_1(\sigma))$ and $a_k$ is above
$h_k$, we have that $d_{t-q_k}$ is above $h_k$.  Notice that the rightmost $1$ of
the obtained $J_{k-1}$ is $h_k$. Thus, this $J_{k-1}$, combining
with $d_{t-q_k}$, will form a $Q_k$ above $a_1$ in $\sigma$, which
contradicts the hypothesis that $\sigma$ contains no $Q_k$ above
$a_1$.

 In the latter case,     $h_{k-2}$ is not above  $a_{k-1}$ since $c_1$ is below $a_k$ (and $b_{k-1}$) and there is no $1's$ inside $E_1$.
If  $h_{k-2}$ is to the left of $a_{k-1}$ (and $b_{k-1}$), then
   by repeating the horizontal slide algorithm, we can obtain a $J_{k-2}$ above $a_1$ in $\sigma$ from the $J_{k-2}$ consisting of
    $h_1, h_2, \ldots, h_{k-2}$. Notice that the rightmost $1$ of the resulting $J_{k-2}$ is below $b_{k-2}$  and to the left of $b_{k-1}$.
     Then, this $J_{k-2}$,
    combining with $b_{k-1}$ and $b_k$, will form an $F_k$ in $\sigma$.  This contradicts the selection of $b_{k-2}$.

Now suppose that  $h_{k-2}$ is either equal to $a_{k-1}$ or is at
the right of $a_{k-1}$ (and $b_{k-1}$), then by the claim obtained
in the proof of Lemma \ref{psi3}, $h_{k-1}$ is above  $b_{k-1}$ and
to the left of $E_1$.  If  $h_{k-2}\neq a_{k-1}$, then $b_1, b_2,
\ldots, b_{k-1}, h_{k-1}$ form an $H_k$ above $a_1$ in $\sigma$,
which contradicts the hypothesis. If  $h_{k-2}=a_{k-1}$, then we
have  $c_1$ is above $b_{k-2}$ (and $a_{k-1}$). Thus,  according to
the definition of $\psi_2$, there must exists $s$ such that $s\in
\mathcal{D}(\sigma)$ and $q_{k-1}\leq s< q_k-1$. Recall that $c_1$
and $h_{k-1}$ lie in columns $q_k$ and $q_{k-1}$, respectively. From
the equality $\mathcal{D}(\sigma)=\mathcal{D}(\psi(\sigma))$, it
follows that $b_1, b_2, \ldots, b_{k-1}, h_{k-1}$ form an $H_k$
above $a_1$ in $\sigma$, which contradicts the hypothesis.
   Hence, we have concluded
that the assertion also holds for $\psi_2(\sigma)$. \qed

 \noindent{\bf The proof of Lemma \ref{psi2}.}   If not, suppose that    $G$ is   an $H_k$  above  $a_1$  in $\psi(\sigma)$. Label the $1's$
 in $G$ by $g_1, g_2, \ldots, g_k$, from left to right.  Moreover, let $H$ be a $Q_k$ above $a_1$ in $\psi(\sigma)$ such that the rightmost two $1's$ lie in two consecutive columns.  We label its $1's$ in $H$ by $h_1, h_2, \ldots, h_k$.  According to the definition of $Q_k$,    there is a  $Q_k$ above $a_1$ in $\psi(\sigma)$ if and only if  there  exists such an $H$.

 We wish to replace some $1's$ of $G$ (resp.  $H$ ) to form an $H_k$ (resp. $Q_k$) above $a_1$ in $\sigma$, which contradicts the hypothesis that there is no $H_k$ (resp. $Q_k$) above $a_1$ in $\sigma$. Here we only consider the case when $G$ (resp. $H$) is in $\psi_2(\sigma)$.
 The other cases can be verified by the similar arguments.  Since the map $\psi_2$  does not change the positions of any other $1's$,  one of $a_i's$ and $c_i's$    must fall in $G$ (resp.  $H$).

 We can     get a $J_{k-2}$ above $a_1$ in $\sigma$ from the $J_{k-2}$ consisting of $g_1, g_2, \ldots, g_{k-2}$ (resp, $h_1, h_2, \ldots, h_{k-2}$),
 by repeating  the horizontal slide algorithm and replacing each $c_j$ by $d_j$ whenever $c_j$ falls in $G$ (resp. $H$).  From Lemmas \ref{psi3} and $\ref{psi4}$,  it follows that the squares   $g_k$ and $g_{k-1}$ (resp. $h_k$ and $h_{k-1}$) are also filled with $1's$ in $\sigma$.
Hence, the obtained $J_{k-2}$ combining with $g_k$ and $g_{k-1}$
(resp. $h_k$ and $h_{k-1}$) will form a $J_k$ (resp, $G_k$) in
$\sigma$.  Since $h_{k-1}$ and $h_k$ lie in two consecutive columns, the obtained $G_k$ is a $Q_k$.  Recall that   $\mathcal{D}(\pi)=\mathcal{D}(\psi(\sigma))$ and the $1's$ positioned at $g_{k-1}$ and $g_k$ belong to two different blocks of $\psi(\sigma)$. This yields that   the $1's$ positioned at $g_{k-1}$ and $g_k$ also belong to two different blocks of $\pi$.  Thus,
the obtained $J_k$  is an $H_k$. This
completes the proof. \qed

\begin{lemma}\label{psi5} If $\sigma$ contains no $H_k$ or $Q_k$ that is above $a_1$, then

\begin{itemize}
\item[(1)] there exists no $1$ that is above and to the left of $a_k$ such that this $1$, combining with $a_1, a_2, \ldots, a_{k-1}$, forms an $H_k$ in $\psi(\sigma)$;
        \item[(2)] there exists no $1$ that is to the left of $a_k$ in $\psi_1(\sigma)$ (or $\psi_4(\sigma)$), such that this $1$,  combining with $a_1, a_2, \ldots, a_{k-1}$, forms a  $Q_k$ in $\psi_1(\sigma)$  ( or $\psi_4(\sigma)$);
\item[(3)] for $1\leq t\leq k-2$,  the board that is above and to the right of $a_{t}$ cannot contain an $H_{k-t}$ or $Q_{k-t}$ in $\psi(\sigma)$ such that the lowest $1$ of this $H_{k-t}$ or $Q_{k-t}$ is to the left of $a_{t+1}$, and this $H_{k-t}$ or $Q_{k-t}$ , combining with $a_1, a_2, \ldots, a_t$, forms an $H_k$ or $Q_k$ in $\psi(\sigma)$.

\end{itemize}
\end{lemma}
\pf (1) If there is such a $1$, then this $1$, combining with $b_1,
b_2, \ldots, b_{k-1}$, forms an $H_k$ in $\sigma$ since
$\mathcal{D}(\sigma)= \mathcal{D}(\psi(\sigma))$. This contradicts
the hypothesis that there is no $H_k$ above $a_1$ in $\sigma$.

(2) The result follows immediately from the fact that there is no
$1's$  below and to the right of $b_{k-1}$ (and $a_{k-1}$), and the
left of $a_k$.

 (3) If not, suppose that $G$ is such an $H_{k-t}$ (or $Q_{k-t}$) in $\psi(T)$. Label its $1's$ by $g_{t+1}, g_{t+2}, \ldots, g_{k}$, from left to right. By hypothesis, $g_{t+1}$ is to the left of $a_{t+1}$. By the same reasoning as in the proof of Lemmas \ref{psi3} and \ref{psi4},  one can verify that both the squares $g_{k}$ and $g_{k-1}$ are also filled with $1's$ in $\sigma$.
 This ensures that by repeating the horizontal slide algorithm and replacing each $c_j$ (resp. $e_j$) by $d_j$ (resp. $f_j$) in
 $\psi_2(\sigma)$ (resp. $\psi_3(\pi)$ and $\psi_4(\pi)$), we can get an $H_{k-t}$ (or $Q_{k-t}$) in $\sigma$,
 in which $g_{t+1}$ is  leftmost $1$. This $H_{k-t}$ (or $Q_{k-t}$), combining with $b_1, b_2, \ldots, b_t$, forms an $H_{k}$ (or $Q_{k}$)
 in $\sigma$, which is above $a_1$. This contradicts the hypothesis that there is no $H_k$ or $Q_k$ above
 $a_1$. This completes the proof. \qed

 \noindent {\bf The proof of Theorem \ref{th3}. }  First, we aim to show that the map  $\Phi$ is well defined, that is, after finitely many iterations of $\phi$, there will be no occurrences of $H_k$ or $Q_k$.   Suppose that we start with some $\tau\in \mathcal{S}_{n}(F_k)$. At the $t$th application of $\phi$ we select a copy of $H_k$ (or $Q_k$) in $\phi^{t-1}(\tau)$.  This has its lowest $1$ in some row $r$. By Lemma \ref{phi3}, the $H_k$ (or $Q_k$) we will select in $\phi^{t}(\tau)$ cannot have its lowest $1$ anywhere  above row $r$. If it is  in row $r$, then we know it is further to the right than at the previous iteration, because there is only one $1$ in that row, and we have just moved it to the right,from $a_1$ to $b_k$. It follows that each iteration the selection of $a_1$ can only go down or slide right, and therefore the map $\Phi$ is well defined.

 Next we aim  to show that $\mathcal{D}(\tau)=\mathcal{D}(\Phi(\tau))$. We prove by induction on $t$.  Suppose that for any $j<t$, we have $\mathcal{D}(\phi^{j-1}(\tau))=\mathcal{D}(\phi^{j}(\tau))$. We wish to show that $\mathcal{D}(\phi^{t-1}(\tau))=\mathcal{D}(\phi^{t}(\tau))$. At the $t$th application of $\phi$ we select a copy of $H_k$ (or $Q_k$) in $\phi^{t-1}(\tau)$.  This has its lowest $1$ in some row $a_1$. Recall that we have shown that each iteration the selection of lowest square of the selected $H_k$ (or $Q_k$) can only go down or slide right.
  By Lemma \ref{phi1}, there is no $F_k$ with at least one square below $a_1$ in $\phi^{t-1}(\tau)$. From Lemma \ref{phid}, it follow that $\mathcal{D}(\phi^{t-1}(\tau))=\mathcal{D}(\phi^{t}(\tau))$.

 Now we proceed to show that the map $\Psi$ is the inverse of the map $\Phi$. To this end, it suffices to show that $\psi(\phi^t(\tau))=\phi^{t-1}(\tau)$. For our convenience,   let $\pi=\phi^{t-1}(\tau)$ and $\sigma=\phi^{t}(\tau)$.  Suppose that at the $t$th application of $\phi$ we select a copy of $H_k$ (or $Q_k$) in $\pi$, in which the $1's$ are positioned in the squares
 $(p_1, q_1),  (p_2, q_2), \ldots,  (p_k, q_k)$, from left to right.   We have four cases.

 Case 1.  The selected  $1's$    form  a copy of $H_k$,  and $\pi_{q_k-1}>\pi_{q_k+1}$ or $q_k=n$.
 In this case, find the largest $s$ such that $q_{k-1}<s<q_k$ and $s-1\in \mathcal{D}(\pi)$.
 By the construction of the transformation  $\phi$,  the squares $(p_2, q_1) (p_3, q_2 ), \ldots,$ $ (p_k, q_{k-1}), (p_1, s), (\pi_{s}, s+1), \ldots, (\pi_{q_k-1}, q_k)$ are filled with $1's$ in $\sigma$,  and all the other rows and columns are the same as $\pi$.   Note that the $1's$ positioned at  the squares $(p_2, q_1) (p_3, q_2 ), \ldots,$ $ (p_k, q_{k-1}), (p_1, s)$  form an $F_k$ in $\sigma$. Lemmas  \ref{phi1}   and \ref{phi2}    ensure that when we apply the map $\psi$ to $\sigma$, the squares we selected are just $(p_2, q_1) (p_3, q_2 ), \ldots,$ $ (p_k, q_{k-1}), (p_1, s)$.
   By Lemma \ref{phi3}, there is no $H_k$ or $Q_k$ above  row $p_1$. This implies that $\psi(\sigma)$ is well defined.
   Suppose that $\sigma=\{(\sigma_1, 1), (\sigma_2, 2), \ldots, (\sigma_n,n)\}$.
   Clearly, we have  $\sigma_{q_i}=p_{i+1}$ for $i=1, 2, \ldots, k-1$,  $\sigma_{s}=p_1$ and $\sigma_j=\pi_{j-1}$ for $j=s+1, s+2, \ldots, q_k$.

We claim that $\sigma_{s-1}>\sigma_{s+1}$. If $s-1\neq q_{k-1}$, then we have $\sigma_{s-1}=\pi_{s-1}$. Since $s-1\in \mathcal{D}(\pi)$,
we have $\pi_{s-1}>\pi_s$. In this case, we have $\sigma_{s-1}=\pi_{s-1}>\pi_s=\sigma_{s+1}$.
If $s-1= q_{k-1}$, then we have $\sigma_{s-1}=p_{k}$. Recall that we have $\pi_s<\pi_{s+1}<\ldots<p_k$.
This implies that $\sigma_{s-1}=p_{k}>\pi_s=\sigma_{s+1}$. Hence, we have concluded that $\sigma_{s-1}>\sigma_{s+1}$.

      We claim that if $\sigma_{q_{k-1}}<\sigma_{q_{k-1}+1}<\ldots<  \sigma_{s-1}> \sigma_{s}$,
      then we have $\sigma_{q_{k-2}}>\sigma_{s+1}$.  If not, since $\sigma_{q_{k-2}}=p_{k-1}=\pi_{q_{k-1}}$ and $s-1\in \mathcal{D}(\pi)$, we have $s-1\neq q_{k-1}$. Then the $1's$ positioned at the squares
        $(p_2, q_2), (p_3, q_3), \ldots, (p_{k-1}, q_{k-1}), (\pi_{s-1}, s-1), (\pi_{s}, s)$ will form a $Q_k$ above row $p_1$ in $\pi$,  which contradicts the selection of $(p_1, q_1)$.  Hence  the claim is proved.

  Then, according to the definition of map $\psi$,  we have $\psi(\sigma)=\psi_2(\sigma)$.
   Since we have  $\mathcal{D}(\pi)=\mathcal{D}(\sigma)$ and $\pi_s<\pi_{s+1}<\ldots <p_k$, we have $\sigma_{s}<\sigma_{s+1}<\ldots<\sigma_{q_k }$.
    Recall that there are no $1's$ inside $E_1$, we have either $q_k=n$ or $q_k\in \mathcal{D}(\pi)$. This yields that we have either $q_k=n$ or $q_k\in \mathcal{D}(\phi(\pi))$   Hence,
   when we apply the the map $\psi_2$ to $\sigma$, $q_k$ is the largest integer $m$ such that $m>s$ and $m\in \mathcal{D}(\sigma)$ or $m=n$.
      Thus, it is easily seen that $\psi_2(\sigma)=\pi$, that is, $\psi(\phi^t(\tau))=\phi^{t-1}(\tau)$.

 Case 2.  The  selected  $1's$     form  a copy of $H_k$  and $\pi_{q_k-1}<\pi_{q_k+1}$.  In this case, find the least  $t$   such that $t> q_k$ and $t\in \mathcal{A}(\pi)$ or $t=n$.
 By the construction of the map $\phi$,    the squares $(p_2, q_1) (p_3, q_2 ), \ldots, $ $(p_k, q_{k-1}), (\pi_{q_k+1}, q_k), (\pi_{q_k+2}, q_k+1), \ldots, (\pi_{t}, t-1)  (p_1, t)$ are filled with $1's$ in $\sigma$,  and all the other rows and columns are the same as $\pi$.  Note that the $1's$ positioned at  the squares $(p_2, q_1) (p_3, q_2 ), \ldots,$ $ (p_k, q_{k-1}), (p_1, t)$  form an $F_k$ in $\sigma$. Lemmas  \ref{phi1}   and \ref{phi2}    ensure that when we apply the map $\psi$ to $\sigma$, the squares we selected are just $(p_2, q_1) (p_3, q_2 ), \ldots,$ $ (p_k, q_{k-1}), (p_1, t)$.
   By Lemma \ref{phi3}, there is no $H_k$ or $Q_k$ above  row $p_1$. This implies that $\psi(\sigma)$ is well defined.
   Clearly, we have  $\sigma_{q_i}=p_{i+1}$ for $i=1, 2, \ldots, k-1$,  $\sigma_{t}=p_1$ and $\sigma_j=\pi_{j+1}$ for $j=q_k, q_k+1, \ldots, t-1$.

Since $t\in \mathcal{A}(\pi)$ or $t=n$, we have $\pi_{t}<\pi_{t+1}$
or $t=n$.   This implies that
$\sigma_{t-1}=\pi_t<\pi_{t+1}=\sigma_{t+1}$ or $t=n$. By Remark 3.1,
we see that there exits an $s$ such that $s-1\in \mathcal{D}(\pi)$
and $q_{k-1}<s<q_k$. This implies that $\pi_{s-1}>\pi_s<\pi_{s+1}$.
Since $\mathcal{D}(\pi)=\mathcal{D}(\sigma)$, we have
$\sigma_{s-1}>\sigma_s<\sigma_{s+1}$ and $p_{k-1}<s<q_k<t$.
  Then, according to the definition of map $\psi$,  we have $\psi(\sigma)=\psi_3(\sigma)$.
   When we apply the the map $\psi_3$ to $\sigma$, since we have $\mathcal{D}(\pi)=\mathcal{D}(\sigma)$ and  $\pi_{q_k-1}<\pi_{q_k}>\pi_{q_k+1}>\ldots> \pi_{t}$,  $q_k$ is the largest integer $m$
  such that $m-1\in \mathcal{A}(\sigma)$ and $q_{k-1}<m\leq t$.  Thus, it is easily seen that $\psi_3(\sigma)=\pi$, that is, $\psi(\phi^t(\tau))=\phi^{t-1}(\tau)$.

Case 3. The  selected  $1's$     form  a copy of $Q_k$ and $q_k\in \mathcal{A}(\pi)$. By the construction of the map $\phi$,    the squares $(p_2, q_1) (p_3, q_2 ), \ldots, $ $(p_k, q_{k-2}),    (p_1, q_k)$ are filled with $1's$ in $\sigma$,  and all the other rows and columns are the same as $\pi$.  Note that the $1's$ positioned at  the squares $(p_2, q_1) (p_3, q_2 ), \ldots,$ $ (p_k, q_{k-2}), (p_1, q_k)$  form an $F_k$ in $\sigma$. Lemmas  \ref{phi1}   and \ref{phi2}    ensure that when we apply the map $\psi$ to $\sigma$, the squares we selected are just $(p_2, q_1) (p_3, q_2 ), \ldots,$ $ (p_k, q_{k-2}), (p_1, q_k)$.
   By Lemma \ref{phi3}, there is no $H_k$ or $Q_k$ above  row $p_1$. This implies that $\psi(\sigma)$ is well defined.
    Clearly, we have  $\sigma_{q_i}=p_{i+1}$ for $i=1, 2, \ldots, k-3$,  $\sigma_{q_{k-2}}=p_k$ and $\sigma_{q_k}=p_1$.

 According to the definition of $Q_k$, we have $\pi_{q_{k-1}}<\pi_{q_{k-1}+1}<\ldots< \pi_{q_k-1}>\pi_{q_k}$. Moreover,
  we have  $\sigma_{j}=\pi_j $ for $j=q_{k-1}, q_{k-1}+1, \ldots, q_k-1$.
   Thus, we have $\sigma_{q_{k-1}}<\sigma_{q_{k-1}+1}<\ldots<\sigma_{q_k-1}>p_1=\sigma_{q_k}$ and
   $\sigma_{q_{k-2}}=p_k=\pi_{q_k}<\pi_{q_k+1} =\sigma_{q_k+1}$.  Then, according to the definition of map $\psi$,
   we have $\psi(\sigma)=\psi_1(\sigma)$. Thus, it is easily seen that $\psi_1(\sigma)=\pi$, that is, $\psi(\phi^t(\tau))=\phi^{t-1}(\tau)$.

 Case 4.  The  selected  $1's$     form  a copy of $Q_k$,  and $\pi_{q_k}>\pi_{q_k+1}$ or $q_k=n$. In this case, let
    $t$ be the least  such that $t> q_k$ and $t\in \mathcal{A}(\pi)$ or $t=n$. By the construction of the map $\phi$,
     the squares $(p_2, q_1) (p_3, q_2 ), \ldots, $ $(p_k, q_{k-2}), (\pi_{q_k+1}, q_k), (\pi_{q_k+2}, q_k+1), \ldots, (\pi_{t}, t-1)  (p_1, t)$
      are filled with $1's$ in $\sigma$,  and all the other rows and columns are the same as $\pi$.
       Note that the $1's$ positioned at  the squares $(p_2, q_1) (p_3, q_2 ), \ldots,$ $ (p_k, q_{k-2}), (p_1, t)$
       form an $F_k$ in $\sigma$. Lemmas  \ref{phi1}   and \ref{phi2}    ensure that when we apply the map $\psi$ to $\sigma$,
        the squares we selected are just $(q_2, p_1) (q_3, p_2 ), \ldots,$ $ (p_k, q_{k-1}), (p_1, t)$.
   By Lemma \ref{phi3}, there is no $H_k$ or $Q_k$ above  row $p_1$. This implies that $\psi(\sigma)$ is well defined.
   Clearly, we have  $\sigma_{q_i}=p_{i+1}$ for $i=1, 2, \ldots, k-3$,   $\sigma_{q_{k-2}}=p_{k}$, $\sigma_{t}=p_1$ and $\sigma_j=\pi_{j+1}$
   for $j=q_k, q_k+1, \ldots, t-1$.

Since $t\in \mathcal{A}(\pi)$ or $t=n$, we have $\pi_{t}<\pi_{t+1}$ or $t=n$. This implies that
 $\sigma_{t-1}=\pi_t<\pi_{t+1}=\sigma_{t+1}$ or $t=n$.  According to the definition of $Q_k$,
 we have $\pi_{q_{k-1}}<\pi_{q_{k-1}+1}<\ldots< \pi_{q_k-1}>\pi_{q_k}$.  Thus,
  we have
$\sigma_{q_{k-1}} <\sigma_{q_{k-1}+1}<\ldots<\sigma_{q_k-1}>\sigma_{q_k}>\sigma_{q_k+1}>\ldots>\sigma_{t}=p_1$.
  Then, according to the definition of map $\psi$,  we have $\psi(\sigma)=\psi_4(\sigma)$.
  Thus, it is easily seen that $\psi_4(\sigma)=\pi$, that is, $\psi(\phi^t(\tau))=\phi^{t-1}(\tau)$.

 So far, we have deduced that $\psi(\phi^t(\tau))=\phi^{t-1}(\tau)$.

Now we proceed to    to show that the map $\Psi$ is well defined, that is, after finitely many iterations of $\psi$,
there will be no occurrences of $F_k$. Suppose that we start with some $\tau\in \mathcal{S}_{n}(H_k, Q_k)$.
At the $t$th application of $\psi$ we select a copy of $F_k$   in $\psi^{t-1}(\tau)$.
 This has its lowest $1$ in some row $r$.    By Lemma \ref{psi1}, the $F_k$ we will select in $\phi^{t}(\tau)$ cannot
 have its lowest $1$ anywhere below row $r$. If it in row $r$, then we know it is further to the left than at the previous iteration,
 because there is only one $1$ in that row, and we have just moved it to the left, from $b_k$ to $a_1$.
 It follows that  at each iteration the selection of $b_k$ can only go up or slide left.
 Moreover, Lemma \ref{psi2} implies that there is no $H_k$ or $Q_k$ above   $b_k$. Therefore,
 after finitely many iterations of $\psi$, there will be no occurrences of $F_k$.

 Next we aim  to show that $\mathcal{D}(\tau)=\mathcal{D}(\Psi(\tau))$. We prove by induction on $t$.  Suppose that for any $j<t$, we have $\mathcal{D}(\psi^{j-1}(\tau))=\mathcal{D}(\psi^{j}(\tau))$. We wish to show that $\mathcal{D}(\psi^{t-1}(\tau))=\mathcal{D}(\psi^{t}(\tau))$. At the $t$th application of $\psi$ we select a copy of $F_k$  in $\psi^{t-1}(\tau)$.  This has its lowest $1$ in some row $b_k$. Recall that we have shown that each iteration the selection of lowest square of the selected $F_k$   can only go up or slide left.
  By Lemma \ref{psi2}, there is no $H_k$ or $Q_k$ above  $b_k$ in $\psi^{t-1}(\tau)$. Hence,  from Lemma \ref{psi0}, it follow that $\mathcal{D}(\psi^{t-1}(\tau))=\mathcal{D}(\psi^{t}(\tau))$.

  By the same reasoning as in the proof of the equality $\psi(\phi^t(\tau))=\phi^{t-1}(\tau)$, we can prove the equality  $\phi(\psi^t(\tau))=\psi^{t-1}(\tau)$ relying on Lemmas  \ref{psi2} and \ref{psi5}, and the equality $\mathcal{D}(\psi^{t-1}(\tau))=\mathcal{D}(\psi^{t}(\tau))$. The details are omitted here. This completes the proof.
  \qed

\noindent{\bf Acknowledgments.}   The   author was supported by the
National Natural Science Foundation of China.


\end{document}